\documentstyle[12pt]{article}
\begin{document}
\author{S.V. Ludkovsky.}
\title{Duality of $\kappa $-normed topological vector 
spaces and their applications. }
\date{1 November 2000.}
\maketitle
\begin{abstract}
\par  In this article a duality of $\kappa $-normed 
topological vector spaces $X$ is defined and investigated,  
where $X$ is over the field ${\bf K}=
\bf R$ or $\bf C$ or non-Archimedean.
For such spaces the analog of the Mackey-Arens theorem is proved.
A conditional $\kappa $-normability of spaces
$L(X)$ of linear topological homeomorphisms of a locally
convex $\kappa $-normed space $X$ is studied, when an image 
of elements under the corresponding operations is in $L(X)$.
There are investigated cases, when $\kappa $-normability 
of a topological vector space implies its local convexity.
There are given applications of $\kappa $-normed spaces
for resolutions of differential equations
and for approximations of functions in mathematical economy.
\end{abstract} 
\section{Introduction.}
\par  New class of $\kappa $-metric topological spaces
$X$ was introduced earlier in works \cite{shep1,shep2}.
$\kappa $-metric spaces may be non-metrizable, where
a regular $\kappa $-metric is defined as non-negative function
$\rho : X\times 2^X_o\to \bf R$ and satisfying axioms 
$(N1-5)$, where $2^X_o$ is the family of all canonical closed 
subsets of $X$ and $2^X_{\delta }$ is the family of all closed 
$G_{\delta }$-subsets of $X$. 
One of the most important examples of 
$\kappa $-metric spaces are locally compact groups and 
generalized loop groups \cite{luumn48}.
Different distance functions for subsets of linear normed spaces
were studied in article \cite{att}.
Topological vector spaces with $\kappa $-metrics satisfying additional 
conditions related with linearity of these spaces
were defined and studied in work \cite{lusmz}. 
For such $\kappa $-normed spaces (see Defintion 2.1) 
were proved analogs of theorems about fixed point, closed graph and 
open mapping. Free $\kappa $-normed spaces generated by
$\kappa $-metric uniform spaces with uniformly continuous  
$\kappa $-metrics also were studied. There were investigated 
categorial properties of 
$\kappa $-normed spaces relative to products, projective and 
inductive limits.
\par On the other hand, dual pairs and the theorems about dual topologies
play very important role in the theory of topological vector spaces
\cite{nari,schae}.
In this work the duality of $\kappa $-normed spaces is investigated.
Certainly the given definiton of duality for $\kappa $ normed spaces 
differs from that of ordinary topological vector spaces
(see Defintion 2.2).
Theorem 2.4 is the development of the Mackey-Arens theorem
on the case of $\kappa $-normed spaces. Theorem 2.6 is devoted
to the conditional $\kappa $-normability
of the uniform space $L(X)$ of linear topological homeomorphisms
of the $\kappa $-normed space $X$.
Since $L(X)$ is not the linear space, hence 
$\kappa $-normability of $L(X)\times 2^{L(X)}_{\delta }$ 
is considered under conditions, when the corresponding 
operations from Definition 2.1 do not lead outside
the space $L(X)\times 2^{L(X)}_{\delta }$. In particular
$L(X)\times 2^{L(X)}_{\delta }$ is the regular $\kappa $-metrizable space.
Topology of $\kappa $-normed spaces is investigated in
\S 4. There are considered cases, when the $\kappa $-normability
of $(X,2^X_o)$ leads to the $\kappa $-normability
of $(X,2^X_{\delta })$. Then it is proved that the $\kappa $-normability
of $(X,2^X_{\delta })$ implies local convexity of $X$.
In \S 5 of this article applications of $\kappa $-normed spaces
for resolutions of differential equations in infinite-dimensional 
over the field ${\bf K}={\bf R}$ or $\bf C$ of complete
non-normed locally convex spaces $X$ are given.
Evidently the $\kappa $-normed spaces can be used for resolutions 
of more complicated differential equations including partial 
differential equations.  
In \S 6 applications of $\kappa $-normed spaces for approximations 
of functions useful in mathematical economy are studied.
\section{Duality of $\kappa$-normed spaces.}
\par  {\bf 2.1. Defintion.} A topological vector space
$X$ over the filed ${\bf K}=\bf R$ or $\bf C$ or non-Archimedean 
supplied with a family either $S_X:=2^X_o$
of all canonical closed subsets or $S_X=2^X_{\delta }$
of all closed $G_{\delta }$-subsets is called $\kappa $-normed, 
if on the following product $X\times S_X$ 
there exists a non-negative function $\rho (x,C)$,
called a $\kappa $-norm satisfying the following conditions: 
\par $(N1)$ $\rho (x,C)=0$ if and only if $x\in C$;
\par $(N2)$ if $C\subset C'$, then $\rho (x,C)\ge \rho (x,C')$;
\par $(N3)$ $\rho (x,C)$ is uniformly continuous 
by $x\in X$ for each fixed $C\in S_X$;
\par $(N4)$ for each increasing transfinite sequence
$\{ C_{\alpha } \} $ with $C:=cl(\bigcup_{\alpha }C_{\alpha })\in S_X$
the following equality  $\rho (x,C)=\inf_{\alpha }\rho (x,C_{\alpha })$
is satisfied, where $cl_X(A)=cl(A)$ denotes the closure
of a subset $A$ in $X$;
\par $(N5)$ $(a)$ $\rho (x+y,cl(C_1+C_2))\le
\rho (x,C_1)+\rho (y,C_2)$ and 
\par $(N5)(b)$ $\rho (x,C_1)\le \rho (x,C_2)
+{\bar \rho }(C_2,C_1)$ 
(or with the maximum instead of the sum on the right sides of 
inequalities in the non-Archimedean case), where 
$\bar \rho (C_2,C_1):=\sup_{x\in C_2} \rho (x,C_1)$; 
\par $(N6)$ $\rho (\lambda x, \lambda C)=
|\lambda |\rho (x,C)$ for each ${\bf K}\ni \lambda \ne 0$;
\par $(N7)$ $\rho (x+y,y+C)=\rho (x,C)$.
If to consider the empty set $\emptyset $ as the lement of $S_X$, 
then 
\par $(N8)$ $\rho (x,\emptyset )=\infty $ and
$\rho (x,C)<\infty $ for each $x\in X$ and $\emptyset \ne 
C\in S_X$. As usually a subset $A$ in $X$ is called 
$G_{\delta }$ if it is a countable intersection of open subsets.
\par The space $(X,S_X):=X\times S_X$ is called $\kappa $-normed, 
if there is given the fixed $\kappa $-norm $\rho $.
We denote the $\kappa $-normed space by $(X,S_X,\rho )$.
\par Let $X$ be a completely regular topological space, 
then the function 
$\rho $ on $X\times 2^X_o$, satisfying conditions $(N1-N4,N5(b))$
with continuity instead of uniform continuity in Axiom
$(N3)$ is called the regular $\kappa $-metric. 
In works \cite{shep1,shep2} these axioms were denumerated by
$(K1-K5)$ correspondingly.
\par From this definition it follows, that each
normed space is also $\kappa $-normed, if to put $\rho (x,C)=\inf_{y\in C}
\| x-y\|$. The condition $C:=cl(\bigcup_{\alpha }C_{\alpha })\in S_X$
in Axiom $(N4)$ for $S_X={2^X}_o$ is satisfied automatically, since
$C=cl(\bigcup_{\alpha }Int(C_{\alpha }))$.
If $\{ 0 \} \in 2^X_{\delta }=S_X$, then taking $\lambda \to 0$ 
in $(N6)$ implies the trivial equality $\rho (0, \{ 0 \} )=0.$
\par {\bf 2.2. Notes and definitions.} Let $X$ and $Y$ be two locally convex 
spaces over the field ${\bf K}=\bf R$ or $\bf C$ or the non-Archimedean 
spherically complete field $\bf K$ with the non-trivial multiplicative 
norm in this field, such that
these spaces form a dual pair $(X,Y)$. This means that a bilinear functional
$<x,y>$ is given such that $X$ separates points in $Y$ 
(that is, from $<x,y>=0$ for each $y\in Y$
it follows $x=0$) and $Y$ separates points in $X$, where 
$x\in X$, $y\in Y$.
We suppose that the Hausdorff topology $\sf T$ in
$X$ is the topology of the dual pair. The latter means that 
$Y$ is the continuous dual space of $(X,{\sf T})$ 
(see Definition 9.6 and \S 9.202 \cite{nari}). 
\par Then we define the $\kappa $-duality,
that is, a mapping 
$<(x,A)|(y,B)>$ from $(X,S_X)\times (Y,S_Y)$ into
$[0,\infty )$, satisfying the following Conditions $(D1-8)$:
\par $(D1)$ $(a)$ the restriction of $<*|*>$ on 
$(X\times M)\times (Y\times M^{\perp })$
is such that $<(x,M)|(y,M^{\perp })>=|<x+M,y+M^{\perp }>|$, 
where $M\in S_X$ and 
$M^{\perp }\in S_Y$ are $\bf K$-linear subspaces, 
$M^{\perp }:=\{ z\in Y: z(M)= 
\{ 0\} \} $, $<x+M,y+N>$ is the bilinear functional on $(X/M)\times (Y/N)$
for each pair of linear subspaces $M\in S_X$ and $N\in S_Y$;
\par $(D1)(b)$ for each $b\in Y/N$ (or $a\in X/M$), 
if $<a,b>=0$ for each $a\in X/M$
(or $b\in Y/N$), then $b=0$ (or $a=0$ respectively);
\par $(D1)(c)$ for each $x\notin A$, $A\in S_X$ (or $y\notin B$, $B\in S_Y$)
there are $y\in Y$ and $B\in S_Y$ (or $x\in X$ and $A\in S_X$ 
correspondingly) such that $<(x,A)|(y,B)> >0$;
\par $(D2)$ $<(x,A_1)|(y,B_1)>\le <(x,A_2)|(y,B_2)>$ for each
$A_2\subset A_1$ and $B_2\subset B_1$;
\par $(D3)$ $<(x,A)|(y,B)>$ is the uniformly continuous function 
by $x\in X$ 
(or $y\in Y$) for each fixed $A\in S_X$, $y\in Y$ and $B\in S_Y$
(or $x\in X$, $A\in S_X$ and $B\in S_Y$ correspondingly);
\par $(D4)$ $\inf_{\alpha } <(x,A_{\alpha })|(y,B)>=<(x,A)|(y,B)>=
\inf_{\beta } <(x,A)|(y,B_{\beta })>$ for each $x\in X$, $y\in Y$ 
and increasing transfinite sequences $\{ A_{\alpha } \}
\subset S_X$ and $\{ B_{\beta }\} \subset S_Y$ with $A=cl 
(\bigcup_{\alpha }A_{\alpha })$, $B=cl (\bigcup_{\beta }B_{\beta })$,
$A\in S_X$, $B\in S_Y$;
\par $(D5)$ $(a)$ $<(x_1+x_2, cl (A_1+A_2))| (y,B)>\le 
<(x_1,A_1)|(y,B)>+<(x_2,A_2)|(y,B)>$ and 
$<(x,A)|(y_1+y_2, cl (B_1+B_2))>\le
<(x,A)|(y_1,B_1)>+<(x,A)|(y_2,B_2)>$;
\par $(D5)(b)$ $<(x,A)|(y,B)>\le <(x,C)|(y,B)>+
\sup_{z\in C}<(z,A)|(y,B)>$ and
$<(x,A)|(y,B)>\le <(x,A)|(y,E)>+\sup_{z\in E}
<(x,A)|(z,B)>$ for each $x, x_1, x_2 \in X$, $y, y_1, y_2\in Y$,
$A, A_1, A_2, C\in S_X$, $B, B_1, B_2, E\in S_Y$
or with the maximum instead of the sum in the right sides 
of inequalities in the non-Archimedean case;
\par $(D6)$ $<(\lambda x,\lambda A)|(y,B)>=|\lambda |<(x,A)|(y,B)>=
<(x,A)|(\lambda y,\lambda B)>$ for each $x\in X$, $A\in S_X$, $y\in Y$,
$B\in S_Y$, $0\ne \lambda \in \bf K$;
\par $(D7)$ $<(x,A)|(y,B)>=<(x+a,A+a)|(y+b,B+b)>$ for each $a$ and $x\in X$,
$b$ and $y\in Y$, $A\in S_X$, $B\in S_Y$;
\par $(D8)$ $<(x,A)|(y,\emptyset )>=\infty $ and
$<(x,\emptyset )|(y,B)>=\infty $ for each $x\in X\setminus A$, $A\in S_X$, 
$y\in Y\setminus B$, $B\in S_Y$ such that $<(x,A)|(y,B)><\infty $, if
$A\ne \emptyset $ and $B\ne \emptyset $.
\par The mapping $<(x,A)|(y,B)>$ we call the $\kappa $-form.
Usually we consider $S_X=2^X_{\delta }$ and $S_y=2^Y_{\delta }$. 
This $\kappa $-form can be defined for $2^X_o$ and $2^Y_o$
instead of $2^X_{\delta }$ and $2^Y_{\delta }$, but 
$(D1)(a)$ may be not satisfied, since there are locally 
convex spaces $X$ and $Y$ when all proper subspaces 
of $X$ and $Y$ may not belong to $2^X_o$ or $2^Y_o$.
\par {\bf 2.3. Lemma.} {\it For each $\kappa $-duality
of $(X,2^X_{\delta })$ with $(Y,2^Y_{\delta })$ 
from the $\kappa $-normability of $(X,2^X_{\delta })$ the
$\kappa $-normability of $(Y,2^Y_{\delta })$ follows.}
\par {\bf Proof.} Let $\rho _X$ be a $\kappa $-norm
in $(X,2^X_{\delta })$ and $M\in 2^X_{\delta }$
be a $\bf K$-linear subspace in $X$,
then in accordance with Lemma 2 from \cite{lusmz} there is the equality
$M=\bigcap_jC_j$, where $C_j\in 2^X_o$, $C_{j+1}
\subset Int (C_j)$ for each $j$. From Axioms $(N1,3,5-7)$ for $\rho _X$ 
it follows, that $\| x+M \|_{X/M}:=\rho _X(x,M)$ is the norm
in the quotient space $X/M$. Terefore $M$ is equal to 
$M=\bigcap_n\lambda _nV$,
where $U$ is a convex balanced neighbourhood of zero in $X$, 
$U=\theta ^{-1}_M
\{ a\in X/M: \| a\| _{X/M}<1 \}$, $\theta _M: X\to X/M$ is the quotient 
mapping, $V=M+U$, $0\ne \lambda _n\in \bf K$, $\lim_n\lambda _n=0$.
\par Locally convex spaces $X$ and $Y$ form the dual pair, 
hence $M^o=M^{\perp }=\bigcap_n
\lambda _nQ^o=(\bigcup_n \lambda _n^{-1}Q)^o$, where
$Q$ is a convex balanced absorbing bounded subset in
$M$, $Q\ni 0$, $Q^o:= \{ b\in Y:$ $|<q,b>|\le 1$ $\mbox{ for each } 
q\in Q \} $ is the (absolute) polar of a subset $Q$.
Consequently, $M^{\perp }\in 2^Y_{\delta }$, since $M^{\perp }
\subset Q^o+M^{\perp }\subset Q^o$, also $(Q^o+M^{\perp })$ contains 
the unit ball $B(Y/M^{\perp },0,1)$ from $Y/M^{\perp }$
and $M^o$ is closed in $Y$, where
\par $\| b\|_{Y/M^{\perp }}:=\sup_{0\ne a\in X/M}
( |<a,b>|/ \| a\|_{X/M})$, \\ 
$B(Y,y,r):= \{ z\in Y: \| z-y \|_Y \le r \} $
(see \S 9.3, \S 9.8 and \S 9.202 \cite{nari}).
Then we put $\rho _Y(y,B):=\sup_{\rho _X(x,A)>0}<(x,A)|(y,B)>/\rho _X(x,A)$ 
for each $y\in Y$ and $B\in 2^Y_{\delta }$. Now we need to verify that
$\rho _Y$ satisfies Axioms $(N1-N8)$.
\par $(N1)$ In view of $(D1c)$ for each $y\notin B\in 2^Y_{\delta }$ 
there are $x\notin A\in 2^X_{\delta }$ with 
$<(x,A)|(y,B)> >0$, consequently,
$\rho _Y(y,B) >0$. If $\rho _Y(y,B)=0$, then $<(x,A)|(y,B)>=0$
for each $(x,A)\in X\times S_X$, hence $y\in B$.
\par $(N2)$ In view of $(D2)$ we get: $\rho _Y(y,B_1)=\sup_{x\notin A}
<(x,A)|(y,B_1)>/\rho _X(x,A)\le \sup_{x\notin A}<(x,A)|(y,B_2)>/
\rho _X(x,A)=\rho _Y(y,B_2)$ for each $B_2\subset B_1\in 2^Y_{\delta }$.
\par The satisfaction of $(N4-8)$ follows from $(D4-8)$.
\par $(N3)$. Due to $(D3)$ for each $\epsilon >0$, $x\in 
X\setminus A$, $A\in 2^X_{\delta }$ and $B\in 2^Y_{\delta }$
there exists a neighbourhood $F$ of zero
in $Y$ such that $|<(x,A)|(y_1,B)>-<(x,A)|(y_2,B)>|/\rho _X(x,A)<\epsilon $
for each $y_1-y_2\in F$. Let $N=\bigcap_n\lambda _nInt( B)$
for $B\in 2^X_o$, where $Int(E)$ is the interior of a subset $E$ in $X$. 
Then $N\in 2^X_{\delta }$ and $N$ is the closed
$\bf K$-linear subspace. Let $E$ be a $\sigma (Y,X)$-compact disk 
(or a $\sigma (Y,X)$-bounded $c$-compact $\bf K$-disk in the 
non-Archimedean case)
in $Y$, then $E\cap N$ is a $\sigma (N,X)$-compact 
disk (or a $\sigma (N,X)$-bounded $c$-compact $\bf K$-disk
correspondingly) in $N$. 
The Mackey-Arens theorem states, that for the dual pair
$(X,Y)$ the locally convex Hausdorff topology $\sf T$ in $X$ 
is the topology of the dual pair if and only if $\sf T$ is 
the polar topology defined with the help of the family $\sf G$ 
of $\sigma (Y,X)$-compact disks 
(or $\sigma (Y,X)$-bounded $c$-compact $\bf K$-disks)
in $Y$, which cover $Y$
(see the Mackey-Arens theorem in \S (9.6.2) and \S 9.202 in \cite{nari}).
Then $(E\cap N)^o\supset E^o$ is the neighbourhood of zero in $X$,
$Int(E^o)\ni 0$.
We take $E$ from the family $\sf G$ such that $\bigcup_{E\in \sf G}E=Y$. 
The Alaoglu-Bourbaki theorem states that
the polar $U^o$ is $\sigma (X',X)$-compact, if $U$ is a neighbourhood of
zero in the topological vector space $X$ (see \S \S (9.3.3) 
and \S 9.202 in \cite{nari}).
Then there exists $E\in \sf G$ for which
$N^{\perp }=N^o=\bigcap_n\lambda _n(E\cap N)^o=(\bigcup_n\lambda _n^{-1}
(E\cap N))^o$ and $Int( E^o)\ni 0$ due to the Alaoglu-Bourbaki theorem
and the equality $E^{oo}=cl _{\sigma (Y,X)}E$
in accordance with \S (9.3.2) and \S 9.202 \cite{nari}.
Consequently, $N^{\perp }\in 2^X_{\delta }$. Then $\rho _X(x,
N^{\perp })$ is the norm in $X/N^{\perp }$ and 
$\rho _Y(y,N)$ is the norm in $Y/N$. Since $B+N=B$, then
due to $(D1)(a,b)$ we  have that $\rho _Y(y,B)$ is uniformly continuous
by $y\in Y$, since due to $(D3)$ and $(N2,6,7)$ there exists a constant $C>0$
such that 
$|<(x,A)|(y_1,B)>-<(x,A)|(y_2,B)>|/\rho _X(x,A)\le C\rho _Y(y_1-y_2,N)$,
if a neighbourhood $F$ is given with the help of the norm
$\rho _Y(y,N)$ in $Y/N$ and the quotient mapping 
$\theta _N: Y\to Y/N$ (see \S 6.2.2 and \S 6.205 in \cite{nari}). 
Lemma 2 \cite{lusmz} states that  for the $\kappa $-mertizable space 
$(X,2^X_{\delta })$ for each $C\in 2^X_{\delta }$ there exists
a sequence $C_j\in 2^X_o$ such that $Int(C_j)\supset C_{j+1}$ 
for each $j\in \bf N$ and $\bigcap_jC_j=\bigcap_jInt(C_j)=C$.
In accordance with the Corollary 5 \cite{lusmz} for the
$\kappa $-normed locally convex space $X$ 
for each $C\in 2^X_{\delta }$ and $\epsilon >0$ there exists
a convex balanced neighbourhood $U\ni 0$ such that
\par $0\le \rho (x,C)-\rho (x,cl(C+U))<\epsilon $ for each $x\in X$.
Vice versa the latter property for the continuous $\kappa $-norm
without Axiom $(N3)$ implies its uniform continuity.
Then if $x\in B\in 2^X_{\delta }$, then $0\in (B-x)\in 2^X_{\delta }$.
With the help of Axiom $(N4)$, Lemma 2 and Corollary 5
from \cite{lusmz} and the equality
$(\bigcup_{n=1}^{\infty }C_n)^o=\bigcap_{n=1}^{\infty }
(C_n^o)$ while $C_n\in 2^X_0$ for each $n$ we get that
$(N3)$ is satisfied on $S_Y=2^Y_{\delta }$.
If $0\in E\in 2^Y_{\delta }$, then for the locally convex space
$Y$ from $E=\bigcap_{n=1}^{\infty }E_n$ with 
$0\in E_n\subset cl(E_n)\in 2^Y_o$
and open $E_n$ it follows, that there are convex balanced
open subsets $V_n$ in $E_n$ for each $n$.
Then $\bigcap_nV_n=:V$ is the convex balanced $G_{\delta }$-subset 
in $Y$, moreover, $V\subset E$, consequently,
the linear subspace $\bigcap_n\lambda _nV=:M$
is the closed $G_{\delta }$-subset in $Y$, since
$M\subset E$, where $\lim_{n\to \infty }\lambda _n=0$, 
$0\ne \lambda _n\in \bf K$ for each $n\in \bf N$.
\par Then we get the following stronger result, where the 
supposition about existence of $\kappa $-duality is omitted.
\par {\bf 2.4. Theorem.} {\it Let $X$ and $Y$ be a dual pair
of locally convex spaces and let $(X,2^X_{\delta })$ be $\kappa $-normable, 
then $(Y,2^Y_{\delta })$ is also $\kappa $-normable.}
\par {\bf Proof.} In view of Lemma 2.3 it is sufficient
to construct a $\kappa $-duality. In view of the Mackey-Arens theorem
the base $\sf B$ of topology $\sf T$ on $X$ is given by polars $E^o$, 
$E\in \sf G$, where
$\sf G$ is a family of $\sigma (Y,X)$-compact diks (or $\sigma 
(Y,X)$-bounded $c$-compact $\bf K$-disks) in $Y$,
for which $\bigcup_{E\in \sf G}E=Y$. In view of the Weak
Representation Theorem (9.2.3)
and \S 9.202 \cite{nari} for each continuous linear functional
$g$ on $X$ there is the unique element $y\in Y$ such that
$g(x)=<x,y>$. Let
\par $(i)$ $<(x,A)|(y,B)>:=\inf_{(a\in A, b\in B)}|<(b-y),(a-x)>|$ for each
$x\in X$, $\emptyset \ne A\in 2^X_{\delta }$, $y\in Y$, 
$\emptyset \ne B\in 2^Y_{\delta }$ and
$<(x,A)|(y,B)>=\infty $ for $A=\emptyset $ or $B=\emptyset $
with $x\notin A$ and $y\notin B$. For the proof it is sufficient
to verify that Conditions $(D1-D8)$ are satisfied.
\par $(D1)(a)$. If $M$ and $N=M^{\perp }$ are $\bf K$-linear 
subspaces in $X$ and $Y$, $M\in 2^X_{\delta }$, $N\in 2^Y_{\delta }$, 
then $<(x,M)|(y,N)>=\inf_{(a\in M,b\in N)}|<(b-y),(a-x)>|=<(x+c,M)|(y+d,N)>$
for each $c\in M$ and $d\in N$, $<(x,M)|(y,N)>=|<\theta _N(y),
\theta _M(x)>|$, where $\theta _M: X\to X/M$ and $\theta _N: Y\to Y/N$ 
are quotient mappings such that
$<x+M,y+N>=<\theta _N(y), \theta _M(x)>$.
\par $(D1)(b)$. Since $Y$ separates points in $X$, and $X$ in $Y$,
then the same is true for the pair $(X/M,Y/N)$, since 
$N^{\perp }=M^{\perp \perp }=M$.
\par $(D1)(c)$. If $E\in 2^X_{\delta }$, then $z+E\in 2^X_{\delta }$ 
for each $z\in X$.
If $x\in X\setminus A$, $A\in 2^X_{\delta }$ is convex and balanced 
(each $E\in 2^X_{\delta }$ with $0\in E$ contains also
such subset $A$ due to local convexity of $X$, if 
$E\in 2^X_0$, we get $A\in 2^X_o$), then
$M=\bigcap_n\lambda _nA$ is the closed $\bf K$-linear subspace,
$x\notin M$, $\theta _M(A)$ is bounded in $X/M$. 
There exists a ball $S$ of radius 
$\epsilon >0$ with the centre in $x+M$ in the normed space
$X/M$ such that $S\cap \theta _M(A)=\emptyset $. 
In view of the Hahn-Banach theorem 
(see \S (8.4.6) and \S 8.203 in \cite{nari})
for the locally convex space $X$ over $\bf R$, $\bf C$ 
or the spherically complete non-Archimedean field
$\bf K$ and the linear subspace $Z$ in $X$ 
each continuous linear functional on $Z$ has a continuous 
extension on $X$. Therefore, there exists a continuous linear 
functional $y\in Y$ such that $y(x)\notin y(A)$, then there exists
a neighbourhood $Int( B) \ni y$ such that $B\in 2^Y_o$
for which $B(x)\cap B(A)=\emptyset $,
since the field ${\bf K}$ is one of the following:
$\bf R$, $\bf C$ or the spherically complete non-Archimedean field,
where $B(A):= \{ b(a): b\in B, a\in A \} $ for $A\subset X$ 
and $B\subset Y$.
Since cokernels of linear functionals are one-dimensional
over $\bf K$, then there are $B\in 2^Y_{\delta }$ 
with $<(x,A)|(0,B)> >0$.
If $x\in A$ or $y\in B$, then $<(x,A)|(y,B)>=0$, since
$0(A-x)=\{ 0\} $, $(B-y)(0)= \{ 0\} $.
\par $(D2)$. $\inf_{(a\in A_1, b\in B_1)}|<(b-y),(a-x)>|\le \inf_{(
a\in A_2,b\in B_2)}|<(b-y),(a-x)>|$ 
for each $A_2\subset A_1\in 2^X_{\delta }$ 
and $B_2\subset B_1\in 2^Y_{\delta }$.
\par $(D4)$. For each $a\in A$ and $b\in B$ there exist
transfinite sequences of elements $a_{\alpha }\in A_{\alpha }$ and
$b_{\beta }\in B_{\beta }$ such that $a_{\alpha }$ converges to $a$,
$b_{\beta }$ to $b$. From the separate continuity of $<(b-y),(a-x)>$ 
by $a\in X$ and also by $b\in Y$ it follows $(D4)$. 
Conditions $(D5-8)$ follow from definitions.
\par $(D3)$. In view of $(D1)$ it is sufficient to consider
the case of the pair $(X/M, Y/M^{\perp })$ of normed spaces
and bounded subsets $A\in 2^{X/M}_{\delta }$ and
$B\in 2^{Y/M^{\perp }}_{\delta }$, that is,
$A\subset B(X/M,0,r_1)$, $B\subset B(
Y/M^{\perp },0,r_2)$, $0<r_1<\infty $, $0<r_2<\infty $.
The topological vector space $X$ over the non-Archimedean field
$\bf K$ with the non-trivial norm or a convex balanced subset
$E$ in $X$ is called $c$-closed, if each $\bf K$-convex balanced
base of the filter in $X$ or in $E$ respectively has a limit point.
In view of the Alaoglu-Bourbaki theorem $(A-x)$ and $(B-y)$ 
are $\sigma (X,Y)$ and $\sigma (Y,X)$-compact (or $\sigma (X,Y)$ 
and $\sigma (Y,X)$-bounded and $c$-compact in the non-Archimedean case) 
respectively (see \S (5.6.2), \S 5.204 \cite{nari} and Theorem 
3.1.2 \cite{eng}).
Since the linear functionals $(a-x)\in (A-x)$ and $(b-y)\in (B-y)$
are continuous in the weak topologies $\sigma (Y,X)$ and $\sigma (X,Y)$ 
on $Y$ and $X$ correspondingly, then the value
$\inf_{(a\in A, b\in B)}|<(b-y),(a-x)>|$ 
of $<(b-y),(a-x)>$ is attained on the corresponding elements
$a_0\in A$ and $b_0\in B$ due to Theorem 3.1.10 \cite{eng}
and in the non-Archimedean case due to
\S 5.204 \cite{nari}. From the inequality $|z(g)|\le \| z\|_{Y/M^{\perp }}
\| g\|_{X/M}$ it follows,  that
$|\inf_{(a\in A, b\in B)}|<(b-y_1),(a-x)>|-\inf_{(a\in A, 
b\in B)}|<(b-y_2),(a-x)>| |
\le \| y_1-y_2\|_{Y/M^{\perp }}r_1$ and
$|\inf_{(a\in A, b\in B)}|<(b-y),(a-x_1)>|-
\inf_{(a\in A, b\in B)}|<(b-y),(a-x_2)>| |
\le \| x_1-x_2\|_{X/M}r_2$.
Moreover, $\rho _X(x+M,M)$ is the continuous norm in $X/M$, that is,
there exists a constant $C>0$, for which $\rho _X(x,M)\le C \| x+M \|_{X/M}$
for each $x\in X$, where $\theta _M(x)=x+M\in X/M$. 
Since $\inf_{(a\in A, b\in B)} |<(b-y),(a-x)>|$
is invariant under the substitution of $x$ on $x+f$ with $f\in M$, 
$y$ and $y+q$ with $q\in M^{\perp }$, then $<(x,A)|(y,B)>$ 
satisfies Conditon $(D3)$.
As at the end of \S 3 the $\kappa $-norm $\rho _Y$ has an extension 
from $(Y,2^Y_o)$ on $(Y,2^Y_{\delta })$.
\par {\bf 2.5. Remark.} 
Vice versa, if  there is given $\rho _Y$ in $(Y,2^Y_{\delta })$ and
$\kappa $-duality by the formula $(i)$, then the $\kappa $-norm
\par $(i)$ ${\tilde \rho }_X(x,A)
:=\sup_{\rho _Y(y,B)>0}<(x,A)|(y,B)>/\rho _Y(y,B)$ is equivalent to
$\rho _X$, that is, by the definition there exist constants
$0<C_1<C_2<\infty $,
for which $C_1{\tilde \rho }_X(x,A)\le \rho _X(x,A)\le C_2
{\tilde \rho }_X(x,A)$ for each $A\in 2^X_{\delta }$ and $x\in X$. 
Indeed, for $\bf K$-linear subspaces $A=M\in 2^X_{\delta }$ and
$B=M^{\perp }\in S_Y$ this is evident from the Hahn-Banach theorem, 
since $\| x+M\|_{X/M}=
\rho _X(x,M)$. On the other hand, from Formulas $2.4.(i)$ and $2.5.(i)$ 
it follows, that
$\tilde \rho _X(x,A)\le \rho _X(x,A)$ for each $A\in 2^X_{\delta }$ 
and $x\in X$. Let $M=\bigcap_n\lambda _nA$ for a convex balanced
subset $A$, where $0\ne \lambda _n\in \bf K$
for each $n\in \bf N$ and $\lim_{n\to \infty }\lambda _n=0$. Then
$\theta _M(A)$ is bounded in $X/M$, consequently, for each 
$\epsilon >0$ there exists $a\in A$ such that 
$| \| a-x+M\|_{X/M} -\rho _X(x,A)|
<\epsilon $. In view of the Hahn-Banach theorem there exists $y\in Y$
such that $\theta _{M^{\perp }}(y)\in Y/M^{\perp }$ and 
$|y(a-x)|=\| a-x+M\|_{X/M}$, 
moreover, $|y(x+M)|\le \| x+M\|_{X/M}$ for each $x\in X$, 
where $\theta _M: X\to X/M$ and $\theta _{M^{\perp }}: Y\to Y/M^{\perp }$
are the quotient mappings. Therefore, $\rho _X$ and $\tilde \rho _X$ 
are equivalent $\kappa $-norms. 
\par Certainly, as in the case of normed spaces this does not mean 
that $X$ is necessarily reflexive.
\par {\bf 2.6.} Let $X$ be a locally convex space
and $L(X)$ be an uniform space of linear topological
homeomorphisms $S: X\to X$ supplied with the base $W(U,V)$ of topology
induced from the space ${\sf O}(X)$ of linear continuous operators from
$X$ into $X$ such that $W(U,V):=\{ S
\in {\sf O}(X): S(U)\subset V \} $, where $U$ and $V\in \sf V$,
$\sf V$ is the base of neighbourhoods of zero in $X$. 
\par We say that $L(X)$ is conditionally $\kappa $-normable 
if and only if there exists a $\kappa $-norm, for which
Axioms $(N1-N8)$ are satisfied, when results of operations
$(A,B)\mapsto A+B$; $(\lambda ,A)\mapsto \lambda A$;
$(S_1+S_2)\mapsto cl(S_1+S_2)$; $(A,S)\mapsto A+S$; 
$(\lambda ,S)\mapsto \lambda S$ and $cl (\bigcup_{\alpha }S_{\alpha })=S$
for increasing transfinite sequences belong to
$L(X)$ and $2^{L(X)}_{\delta }$ correspondingly, where $A$ and $B\in L(X)$, 
$S$ and $S_{\alpha }\in 2^{L(X)}_{\delta }$, $0\ne \lambda \in \bf K$.
\par {\bf Theorem.} {\it From the $\kappa $-normability
of $(X,2^X_{\delta })$ the conditional $\kappa $-normability
of $(L(X),2^{L(X)}_{\delta })$ follows, when
$X$ is over ${\bf K}=\bf R$ or $\bf C$ or the non-Archimedean
spherically complete field with the non-trivial norm.}
\par {\bf Proof.} Each locally convex space $X$ is isomorphic with
the projective limit of normed spaces $Y_{\alpha }$, 
that is, $X=pr-\lim \{ Y_{\alpha }, 
\pi ^{\alpha }_{\beta }, \Lambda \} $, where $\pi ^{\alpha }_{\beta }:
Y_{\alpha }\to Y_{\beta }$ are continuous linear mappings from $Y_{\alpha }$
onto $Y_{\beta }$ for each $\alpha >\beta \in \Lambda $, $\Lambda $  
is a directed set, $\pi ^{\alpha }_{\beta }\circ
\pi ^{\beta }_{\gamma }=\pi ^{\alpha }_{\gamma }$ for each
$\alpha \ge \beta \ge \gamma \in \Lambda $, $\pi _{\alpha }: 
X\to Y_{\alpha }$ are quotient mappings (see \S 2.5 \cite{eng},
\S 6.5, Theorem (6.7.2), \S 6.205 \cite{nari}).
Let $p_{\alpha }$ be a seminorm in $X$ such that 
$X_{\alpha }=p^{-1}(0)$ and $Y_{\alpha }=X/X_{\alpha }$. The seminorm 
$p_{\alpha }$ induces the norm ${\hat p}_{\alpha }$ in $Y_{\alpha }$,
then $X_{\alpha }= \bigcap_n\lambda _n Int(U_{\alpha })$, 
where $U_{\alpha }=
\pi _{\alpha }^{-1}(B(X_{\alpha },0,1))$, $0\ne \lambda _n \in \bf K$
for each $n\in \bf N$ and $\lim_n\lambda _n=0$,
$B(Y_{\alpha },0,1)$ is the ball of radius $1$ with center $0$
in $Y_{\alpha }$. Therefore, $X_{\alpha }\in 2^X_{\delta }$ for each
$\alpha \in \Lambda $. If $M$ is a closed 
$\bf K$-linear subspace in $X$, then $A(M)$ is linearly topologically 
isomorphic with $M$ for each $A\in L(X)$. The operator $A$
generates the opearator $\dot A: X/M\to X/M$ such that $\dot A\in L(X/M)$
and $\dot A(x+M)=A(x)+A(M)$. At the same time $X/M$ is the normed space
with the norm $\rho _X(x,M)$, where $x+M=\theta _M(x)$
and $\theta _M: X\to X/M$ is a quotient mapping.
Then the topolgy in $L(X/M)$ is hereditary from the normed space
${\sf O}(X/M)$.
Evidently ${\bf K}^*L(X)=L(X)$, where ${\bf K}^*:={\bf K}\setminus 
\{ 0\} $. Let $\rho _L(A,S):=\sup_{(x\in X\setminus E, E\in 
2^X_o)}\inf_{(B\in S, a\in E)} (\| (\dot B-\dot A)(x-a+M)\| _{X/M}/
\rho _X(x,E))$ for each $S\in 2^{L(X)}_o$, $S\ne \emptyset $, where
$M=\bigcap_n\lambda _nE$, $\dot A$ and $\dot B$ are operaators on
$X/M$, a subset
$E$ is convex and balanced, $0\ne \lambda _n\in \bf K$
for each $n\in \bf N$ such that  
$\lim_{n\to \infty } \lambda _n=0$.
If $S=\emptyset $ we put $\rho _L(A,\emptyset )=\infty $.
It remains to verify Conditions $(N1-8)$, when the results
of operations
$(A,B)\mapsto A+B$; $(\lambda ,A)\mapsto \lambda A$;
$(S_1+S_2)\mapsto cl(S_1+S_2)$; $(A,S)\mapsto A+S$; 
$(\lambda ,S)\mapsto \lambda S$ and $cl (\bigcup_{\alpha }S_{\alpha })=S$ 
for increasing transfinite sequences belong to
$L(X)$ and $2^{L(X)}_{\delta }$ respectively, where $A$ and $B\in L(X)$, 
$S$ and $S_{\alpha }\in 2^{L(X)}_{\delta }$, $0\ne \lambda \in \bf K$.
\par $(N1)$. If $A\in S$, then from $0(x-E)=\{ 0\}$
it follows that $\rho _L(A,S)=0$. If $A\notin S$, then there exists
a convex balanced neigbourhood $W(U,V)$ of zero in ${\sf O}(X)$
such that $(A+W(U,V))\cap S=\emptyset $, where $U=\bigcap_{j=1}^n\pi _{
\alpha _j}^{-1}(B(Y_{\alpha _j},0,r_j))$, $n\in \bf N$, $0<r_j<\infty $,
$V=\bigcap_{i=1}^m\pi _{\beta _i}^{-1}(B(Y_{\beta _i},0,R_i))$, $m\in \bf N$,
$0<R_i<\infty $.
In view of the Hahn-Banach theorem there are $f\in X'$ and $g\in X$
such that
$q_{A,S}(f,g):=\inf_{B\in S}|f(B-A)g|>0$, since $W(U,V)F$  is convex 
in $X$ for each convex neighbourhood $F$ of zero in $X$. Then in view
of continuity of $q_{A,S}(f,g)$ by $f$ and $g$ there exists a neighbourhood
$U$ such that $\inf_{h\in U}q_{A,S}(f,h+g)>0$, consequently, 
$\inf_{(B\in S,a\in E)}
\| (\dot B-\dot A)(a+M)\|_{X/M}>0$, where $E=g+U$.
\par $(N2)$. Follows from the inequality $\inf_{(B\in S_1, a\in E_1)} 
\| (\dot B-\dot A)(x-a+M) \|_{X/M}\le \inf_{(B\in S_2, a\in E_2)} 
\| (\dot B-\dot A)(x-a+M) \|_{X/M}$ for each $S_2\subset S_1$.
\par $(N3)$. Since $X/M$ and ${\sf O}(X/M)$ are normed spaces, 
then $\inf_{(B\in S, a\in E)} \| (\dot B-\dot A)
(x-a+M) \|_{X/M}$ is uniformly continuous by $\dot A\in L(X/M)$ for
fixed $S\in 2^{L(X)}_o$, $x\in X$ and $E\in 2^X_o$.
Moreover, there exists $C=const >0$ such that 
$\rho _X (x,M)\le C \| x+M\|_{X/M}$
for each $x\in X$. Consequently, for each $\epsilon >0$, 
$x\in X\setminus E$, $E\in 2^X_o$, $S\in 2^{L(X)}_o$
there exists a neighbourhood of zero $F$ in ${\sf O}(X)$
such that $|\inf_{(
B\in S, a\in E)} \| (\dot B-\dot A_1)(x-a+M)\|_{X/M}- \inf_{(B\in S, a\in E)} 
\| (\dot B-\dot A_2)(x-a+M)\|_{X/M}|/\rho _X(x,E)<\epsilon $
for each $A_1-A_2\in F$.
\par $(N4)$. For each $s\in S=cl (\bigcup_{\alpha }S_{\alpha })$
there exists a transfinite sequence $s_{\alpha }\in S_{\alpha }$,
converging to $s\in S$, where $S$ and $S_{\alpha }\in 2^{L(X)}_o$ 
for each $\alpha $, $\{ S_{\alpha }: \alpha \} $ is the increasing
transfinite sequence. Then for each
$x\in X$ and $E\in 2^X_o$ there exists $\lim_{\alpha }(\inf_{a\in E}
\| (s_{\alpha }-\dot A)(x-a+M)\|_{X/M})=\inf_{a\in E} \|(s-
\dot A)(x-a+M)\|_{X/M}$, consequently, $\lim_{\alpha }(\inf_{b\in 
S_{\alpha }, a\in E)} \| (\dot B-\dot A)(x-a+M)\|_{X/M})=\inf_{(B\in S, 
a\in E)} \| (\dot B-\dot A)(x-a+M)\|_{X/M}$, since $| \| \dot B_1(y+
M)\|_{X/M}-\| \dot B_2(y+M)\|_{X/M}|\le \| \dot B_1 - \dot B_2\|_{
{\sf O}(X/M)}\| y+M\|_{X/M}$.
\par $(N5)$. For each $x\in X$, $A$, $B$ and $C\in {\sf O}(X)$
we have
$\| (\dot A-\dot B)(x+M)\|_{X/M}\le \|(\dot A-\dot C)(x+M)\|_{X/M}+
\| (\dot B-\dot C)(x+M)\|_{X/M}$, hence $\inf_{(B\in cl(S_1+S_2), a\in E)}
\| (\dot B-\dot A_1-\dot A_2)(x-a+M)\|_{X/M}\le \inf_{(B_1\in S_1, a\in E)}
\| (\dot B_1-\dot A_1)(x-a+M)\|_{X/M}+\inf_{(B_2\in S_2, a\in E)}\| 
(\dot B_2-\dot A_2)(x-a+M)\|_{X/M}$ and inevitably
$\inf_{(B\in S_1, a\in E)}\| (\dot B-\dot A)(x-a+M)\|_{X/M}\le
\inf_{(C\in S_2, a\in E)}\| (\dot C-\dot A)(x-a+M)\|_{X/M}+
\sup_{C\in S_2}\inf_{(B\in S_1, a\in E)}\| (\dot B-\dot C)(x-a+M)\|_{X/M}$
(with the maximum instead of the sum in the right parts of the inequalities 
in the non-Archimedean case).
Taking $\sup_{(x\in X\setminus E, E\in 2^X_o)}$ of both parts
of inequalities divided by $\rho _X(x,E)$ we get $(N5)$.
\par $(N6-8)$ follow from the definition of $\rho _L$. 
As at the end of \S 3 the $\kappa $-norm has the continuous extension
from $L(X)\times 2^{L(X)}_o$ on $L(X)\times 2^{L(X)}_{\delta }$.
\section{Features of a $\kappa $-normed space topology.}
\par {\bf 3.1. Remark.} Let $X$ be a locally convex space
such that $(X,2^X_{\delta })$ is $\kappa $-normable
by a $\kappa $-norm $\rho $. The topology of $X$ is defined by
a family of seminorms $p_{\alpha }$. Then $X_{\alpha }:=p^{-1}_{
\alpha }(0)$ are closed $\bf K$-linear subspaces in $X$.
In view of the $\kappa $-normability we have as above
$X_{\alpha }\in 2^X_{\delta }$ such that $\rho _{\alpha }
(x):=\rho (x,X_{\alpha })$ are seminorms in $X$. 
Each quotient space $X/X_{\alpha }$ is normable by
a norm ${\hat p}_{\alpha }([x]):=\inf_{z\in X_{\alpha }}
p_{\alpha }(x+z)$, where $x\in [x]:=x+X_{\alpha }$.
Due to uniform continuity of $\rho (x,C)$ by $x\in X$ for each fixed
$C\in 2^X_{\delta }$ we have that the family
$ \{ \rho _{\alpha } \} $ defines topology of $X$
and $\rho _{\alpha }$ are seminorms in $X$. In view of $(N1)$ we get that 
${\hat \rho }_{\alpha }([x]):=\rho _{\alpha }(x,X_{\alpha })$ is the norm
in the quotient space $X/X_{\alpha }$ for each $\alpha $, where
$x\in [x]:=x+X_{\alpha }$.
Therefore, we can choose a family 
$$(i) \mbox{ } \{ p_{\alpha }: \alpha \in \Omega \} $$ 
defining topology of $X$ such that
$$(ii)\mbox{ }p_{\alpha }(x)\le {\hat \rho }_{\alpha }(x)$$ 
for each $\alpha \in \Omega $ and $x\in X$.
\par {\bf Theorem.} {\it Let $(X,2^X_{\delta },\rho )$ be the 
$\kappa $-normed space, whose topology is defined by Family 
$(i)$ of seminorms satisfying Condition $(ii)$, where
$$(iii)\mbox{ } \sup_{\alpha \in \Omega }(x,X_{\alpha })<\infty $$ 
for each $\alpha \in \Omega $.
Then for each countable subfamily $ \{ p_{\alpha }: 
\alpha \in \Lambda \} $, where $\Lambda \subset \Omega $
and $card (\Lambda )\le \aleph _0,$
the following function $q_{\Lambda }(x):=\sup_{\alpha }
p_{\alpha }(x)$ is a continuous seminorm in $X$.}
\par {\bf Proof.} In view of Conditions $(ii,iii)$ 
we have $q_{\Lambda }(x)<\infty $ for each $x\in X$.
Each $p_{\alpha }$ is non-negative, hence
$0\le q_{\Lambda }(x)$ for each $x\in X$.
Since each $p_{\alpha }(x)$ is a seminorm, then
$q_{\Lambda }(x+y)=\sup_{\alpha \in \Lambda }p_{\alpha }(x+y)
\le  sup_{\alpha \in \Lambda }p_{\alpha }(x)+
\sup_{\beta \in \Lambda }p_{\beta }(x)=q_{\Lambda }(x)+q_{\Lambda }(y)$
for each $x, y \in X$, which is the triangle inequality 
for $q_{\Lambda }$. In the non-Archimedeean case we have
the strong triangle inequality 
$q_{\Lambda }(x+y)\le \max (q_{\Lambda }(x), q_{\Lambda }(y))$.
Suppose that for some countable subset $\Lambda $ the seminorm
$q_{\Lambda }$ is not continuous in $X$. We consider 
the subspace $Y:=\bigcap_{\alpha \in \Lambda }X_{\alpha },$
then $Y\in 2^X_{\delta }$, since $Y=\bigcap_{j=1}^{\infty }(
\bigcap_{\alpha \in \Lambda }\lambda _jp_{\alpha }^{-1}([0,1]))$, 
where $\lambda _j\ne 0$ for each $j\in \bf N$ and $\lim_j\lambda _j=0$.
If $q_{\Lambda }$ is not continuous, then
$Int_Xq_{\Lambda }^{-1}([0,1])=\emptyset $, since $q_{\Lambda }^{-1}([0,1])$
is absolutely convex in $X$, consequently, there 
exists a vector $0\ne v\in X\setminus Y$ such that
$\rho (v,Y)>0$. On the other hand, $\rho (x,Y)$ is uniformly continous 
by $x\in X$, consequently, there exists a seminorm $p$ such that
$\rho (x,Y)=p(x)$ for each $x\in X$. 
But in view of Conditions $(i-iii)$ and Axiom $(N8)$ we have
$q_{\Lambda }(x)\le \sup_{\alpha \in \Lambda }\rho (x,X_{\alpha })
\le \rho (x,Y)<\infty $ for each $x\in X$. This inequality contradicts 
our supposition that $q_{\Lambda }$ is discontinuous, consequently,
$q_{\Lambda }$ is continuous.
\par {\bf 3.2. Theorem.} {\it Let $X$ be a locally convex space
and let $(X,2^X_o,\rho )$ be a $\kappa $-normed space.
Suppose that a topology of $X$ is defined by a family of seminorms
$\{ p_{\alpha }: \alpha \in \Omega  \} $ such that
\par $(i)$ $q_{\Lambda }(x)<\infty $ for each $x\in X$ and 
$q_{\Lambda }(x)$ is continuous by $x\in X$
for each countable subset $\Lambda \subset \Omega $, where
$q_{\Lambda }(x):=\sup_{\alpha \in \Lambda }p_{\alpha }(x)$.
Then $\rho $ has the extension on $X\times 2^X_{\delta }$.}
\par {\bf Proof.} Let $X_{\alpha }:=p_{\alpha }^{-1}(0)$, then
we have a $\kappa $-norm $\rho _{\alpha }(x,cl(C+X_{\alpha }))
=: {\hat \rho }_{\alpha } ([x],[C])$ in $(Y_{\alpha },2^{Y_{\alpha }}_o)$
for each $\alpha $, where $Y_{\alpha }:=X/X_{\alpha }$
is the quotient space, $[x]=x+X_{\alpha }$, $[C]:=cl(C+X_{\alpha })$.
Each $Y_{\alpha }$ is the normed space with the norm
${\hat p}_{\alpha }([x]):=\inf_{(z\in X_{\alpha })}
p_{\alpha }(x+z),$ where $x\in [x].$
To finish the proof we need the following lemma.
\par {\bf 3.3. Lemma.} {\it Let $Y$ be a normed space and $\rho $ be  a
$\kappa $-norm on $Y\times 2^Y_o$, then $\rho $ has an extension 
on $Y\times 2^Y_{\delta }$.}
\par {\bf Proof.} Let $p$ be a norm in $Y$, then
$p^{-1}([0,1/n])$ is a canonical closed neighbourhood of zero in $Y$
such that $ \{ 0 \} =\bigcap_{n=1}^{\infty }p^{-1}([0,1/n])
\in 2^Y_{\delta }.$ Put $\rho (x, \{ 0 \} ):=\sup_{n\in \bf N}
\rho (x,p^{-1}([0,1/n])).$ For each $n\in \bf N$ according to Axiom $(N2)$
we have the following inequality: 
$\rho (x,p^{-1}([0,1/(n+1)]))\ge \rho (x,p^{-1}([0,1/n])).$
In view of Axiom $(N5)$ we have another inequality:
$\rho (x,p^{-1}([0,1/(n+k)])) - \rho (x,p^{-1}([0,1/n]))
\le {\bar \rho }(p^{-1}([0,1/n]),p^{-1}([0,1/(n+k)])),$
where ${\bar \rho }(A,B):=\sup_{a\in A}\rho (a,B)$ for each 
$A, B\in 2^Y_0$ and $k\in \bf N$. In view of Corollary 5
\cite{lusmz} for each $\epsilon >0$ there exists
$n\in \bf N$ such that ${\bar \rho }(p^{-1}([0,1/n]),p^{-1}([0,1/(n+k)]))
<\epsilon $ for each $k\in \bf N$. Therefore $\rho (x, \{0 \} )$
is the continuous norm in $Y$. Each norm in $Y$ induces a $\kappa $-norm
in $(Y, 2^Y_{\delta })$, consequently, $\rho $ has an extension
on $Y\times 2^Y_{\delta }$.
\par Continuation of the {\bf proof} of Theorem 3.2.
In view of Lemma 3.3 we have an extension ${\hat \rho }_{\alpha }$
on $Y_{\alpha }\times 2^{Y_{\alpha }}_{\delta }$.
Due to Proposition 13 \cite{lusmz} there exists a locally convex space $Z$ 
and a $\kappa $-norm $\eta $ in $(Z,2^Z_{\delta })$ such that 
each $Y_{\alpha }$ is the quotient space of $Z$ and each
$\kappa $-norm $\eta _{Y_{\alpha }}$ in $(Y_{\alpha },
2^{Y_{\alpha }}_{\delta })$ is equivalent to ${\hat \rho }_{\alpha }$.
From the construction of $Z$ (see Conditions $(a-d)$ in the proof of 
Proposition 13 \cite{lusmz}) and Condition $(i)$ of Theorem 3.2
we have that $X$ is the quotient space of $Z$ as the locally convex space,
since each ${\hat \rho }_{\alpha }(x,C)$ is uniformly continuous by
$x\in Y_{\alpha }$ for each fixed $C\in 2^{Y_{\alpha }}_{\delta }$
relative to the norm ${\hat p}_{\alpha }$.
In view of Theorem 6 \cite{lusmz} the $\kappa $-norm $\eta $ in $Z$
induces a $\kappa $-norm ${\tilde \rho }$ in $(X,2^X_{\delta })$.
The restrictions of these two $\kappa $-norms on $X\times 2^X_o$ are 
equivalent and ${\hat {\tilde \rho }}_{\alpha }(x,C)=
{\hat \rho }_{\alpha }$ for each $x\in Y_{\alpha }$ and 
$C\in 2^{Y_{\alpha }}_o$.
\par {\bf 3.4. Remark.} In the case of $X$ equal
to the countable product
of normed spaces Condition $3.2.(i)$ produces the box topology. 
\par {\bf 3.5. Corollary.} {\it Let $X$ be a complete locally convex space
satisfying Condition $3.2.(i)$ with a $\kappa $-norm
$\rho $ in $(X,2^X_o)$, which induces a metric $D$ in $2^X_o$
in accordance with Definition 7 \cite{lusmz}.
Then the completion of $2^X_o$ is equal to $2^X_{\delta }.$}
\par {\bf Proof.} In view of Theorem 3.2 $\rho $ has the extension
on $X\times 2^X_{\delta }$. In accordance with Theorem 9
\cite{lusmz} $(2^X_{\delta },D)$ is complete, where $D$ has the 
natural extension on $(2^X_{\delta })^2$. In view of Corollary 5
\cite{lusmz} we have that $2^X_o$ is dense in $2^X_{\delta }$.
\par {\bf 3.6. Theorem.} {\it Let $(X,2^X_{\delta },\rho )$ 
be a $\kappa $-normed topological vector space.
Then $X$ is locally convex.}
\par {\bf Proof.} Let us suppose that $X$ is not locally convex and 
we consider a base ${\sf B}_o$ of neighbourhoods of zero in $X$.
Without loss of generality it can be supposed that each $C\in {\sf B}_o$ 
is balanced.
Since $X$ is not locally convex there exists $C\in {\sf B}_o$ such that
$\bigcap_{\lambda _j}\lambda _jC=E$ is not the $\bf K$-linear
subspace of $X$, where $\lambda _j\ne 0$ for each $j\in \bf N$ and
$\lim_j\lambda _j=0$, where $X$ is over the field $\bf K$. 
If $x\in  E$, then
$\lambda _jx\in E$ for each $j$, consequently, ${\bf K}x\subset E$.
Indeed, if $E$ is the $\bf K$-linear subspace, then $E$ is closed in $X$, 
since from $cl(\alpha E)\subset \beta E$ for each $0<|\alpha | < |\beta |$
it follows, that for $E$ it can be chosen $\lambda _j\ne 0$
such that $cl(\lambda _{j+1}C)\subset \lambda _jC$
for each $j\in \bf N$. Therefore, $X/E$ would be the $\bf K$-linear 
topological space. In view of Theorem 6 \cite{lusmz} the quotient space
$X/E$ is $\kappa $-normed with the $\kappa $-norm $\rho _{X/E}$.
But $C+E=C$ in $X$ and $\pi _E(C+E)$ is bounded in $X/E$,
where $\pi _E: X\to X/E$ is the quotient mapping.
Since $E\in 2^X_{\delta }$, then $\rho _{X/E}([x],\{ 0 \} )$
is the norm in $X/E$, where $x\in X$, $[x]=x+E=\pi _E(x)\in X/E$,
$[E]= \{ 0 \} \in X/E$. So it would mean that $X/E$
is the normed space, hence $X$ would be a projective limit
of normed spaces, consequently, $X$ is locally convex.
This contardicts our assumption.
\par This means that there exists $C\in {\sf B}_o$ or $C\in 2^X_{\delta }$
for which $E$ is not the linear subspace of $X$.
Therefore, there exists a system of vectors $\{ x_j: j\in \alpha \} $
in $E$ such that $sp_{\bf K} \{ x_j: j\in \beta \} $
is not contained in $E$ for each subset $\beta $ in $\alpha $,
where $\{ x_j: j\in \alpha \} $ is a system of $\bf K$-linearly 
independent vectors in $X$. If ${\bf K}=\bf R$ or ${\bf K}=\bf C$ we take
${\bf L}:={\bf K}$. We can consider the completion of $X$, 
which is also $\kappa $-normed due to uniform continuity of $\rho $.
So we can suppose that $X$ is complete.
Then $\bf K$ has to be uniformly complete for complete $X$.
If $\bf K$ is the non-Archimedean non-locally compact 
field, then it has a locally compact subfield $\bf L$ such that 
$X$ can be considered as the $\bf L$-linear space $X_{\bf L}$.
Spaces $X$ and $X_{\bf L}$ are uniformly homeomorphic and $\bf L$-linearly,
but not $\bf K$-linearly for ${\bf K}\ne \bf L$.
In view of the Hahn-Banach theorem \cite{nari} there
are $\bf L$-linear continuous functionals $e_j$ such that
$e_j(x_j)=1$. Then $B_j:=\{ y: y\in X_{\bf L}, |e_j(y)|_{\bf L}\le 1 \} $
are canonical closed subsets in $X_{\bf L}$. If 
$\alpha $ is countable, then $C_1:=C\cap (\bigcap_{j\in \alpha }B_j)\in 
2^X_{\delta }$ and $E_1:=\bigcap_{j\in \alpha }\lambda _jC_1$ 
is the closed $\bf L$-linear subspace of $X_{\bf L}.$ 
We get the inverse mapping system $ \{ X/E_1, \pi ^{{E'}_1}_{E_1}, 
\Omega \} $, where ${E'}_1>E_1$ if and only if ${E'}_1\subset E_1
\in \Omega $ and ${E'}_1\ne E_1$, $\Omega $ is the family of such $E_1$,
$\pi ^{{E'}_1}_{E_1}: X/{E'}_1\to X/E_1$ are $\bf L$-linear.
Therefore,
$X_{\bf L}=pr-\lim \{ X/E_1 \} $ is locally convex. On the other hand,
$\bf K$ is a the locally convex space over $\bf L$, consequently,
$X$ is locally convex.
\par It remains to consider the case, when there exists 
$C\in {\sf B}_o$ or $C\in 2^X_{\delta }$ such that
for it the corresponding $E$ has $card (\alpha )>\aleph _0$.
In view of Corollary 5 \cite{lusmz} for each $\epsilon >0$
there exists $U\in {\sf B}_o$ such that
$0\le \rho (x,C)-\rho (x,cl(C+U))<\epsilon $ for each $x\in X$.
Let $z\in cl( sp_{\bf K}E)\setminus E\ne \emptyset $,
then $\rho (z,E)>0$. Suppose $z=\sum_{j\in \alpha }\lambda _jx_j$ and 
$x_j\in U$ for each $\lambda _j\ne 0$, where $\alpha $ is an ordinal
due to the Kuratowski-Zorn lemma. 
Certainly this sum expressing $z$ has only countable 
subset $j\in \beta \subset \alpha $ of $\lambda _j\ne 0$.
Due to Lemma 3.3 for each normed space $Z$ its $\kappa $-norm 
has extension to the norm in $Z$, since 
$\{ z \} \in 2^Z_{\delta }$ for each $z\in Z$.
If for each $\epsilon >0$ there 
exists such $U$, then from the equality
$\rho (\lambda y,\lambda cl(E+U))=|\lambda |\rho (y,cl(E+U))$ 
for each $0\ne \lambda \in \bf K$ it follows that 
$0\le \rho (z,E)-\rho (z,cl(E+U))<\epsilon $, consequently,
$\rho (z,E)<\epsilon $, since $\lim_j\lambda _j=0$ and $cl(E+U)$ contains 
$\sum_{\beta \ni j\ge \alpha _0}\lambda _jx_j$ 
for some $\alpha _0 \subset \alpha $
of the cardinality $card (\alpha _0)<\aleph _0$ and 
each ${\bf K}^{\alpha _0}$ is the normed space.
Since $\epsilon $ is arbitrary, then $z\in E$
and inevitably $cl(sp_{\bf K}E)=E$, contradicting our assumption.
\par Therefore, for each $\epsilon >0$ there exists $U\in {\sf B}_o$
such that $x_j\notin U$ for each $j\in \alpha $ and $\bigcap_j \lambda _jU$
is the $\bf K$-linear subspace in $X$. Then such $U$ would give
$C\cap U$ with $\bigcap_jcl( \lambda _j(C\cap U))=:E_0$
and $X=pr-\lim \{ X/E_0 \} $ again contradicting
assumption on $X$, since $E_0$ are closed $\bf K$-linear subspaces in $X$
and $X/E_0$ are normed spaces. 
Indeed, in view of continuity of addition in $X$ for each
$V\in {\sf B}_o$  there are $C\in {\sf B}_o$ and $U\in {\sf B}_o$
such that $cl(C+U)\subset V$. Therefore, $X$ is the 
$\bf K$-locally convex space.
\par {\bf 3.7. Note.} The latter theorem shows, that  there are
topological vector spaces $X$, for which $(X,2^X_{\delta })$ 
are not $\kappa $-normable, for example, non locally convex 
metrizable $X$.
\section{Embeddings of $2^X_o$ and $2^X_{\delta }$.} 
\par {\bf 4.1.} For a complete locally convex space
$X$ let $S_X$ be either a family $2^X_o$ of all
canonical closed subsets or a family $2^X_{\delta }$ of all closed
$G_{\delta }$-subsets. As in \S 3.6 we consider the subfield
${\bf L}$ of the field $\bf K$ and the corresponding
topological $\bf L$-linear space $X_{\bf L}$. Let  ${\sf P}$ be a subset
of the topological dual space $X^*_{\bf L}$ separating 
points of $X_{\bf L}$.
Let $\mu $ be a non-negative Haar measure on $\bf L$.
Two elements $A, B\in 2^X_{\delta }$ we call
$\Xi $-equivalent if and only if $\mu (\pi (A\bigtriangleup B))=0$
for each $\pi \in \sf P$, where as usually $(A\bigtriangleup B)=
(A\setminus B)\cup (B\setminus A)$. This equvalence is written as 
$A\Xi B$.
\par If $X$ is over $\bf R$ or $\bf C$ we take
$\mu $ such that $\mu ([0,1])=1$. 
In the non-Archimedean case there are two variants:
the characteristic of the field $\bf K$ is
$char ({\bf K})=s>0$; $char ({\bf K})=0$, 
then ${\bf K}$ is the extension of the field $\bf Q_s$ 
for the corresponding prime number $s$, then we take $s\ne p$.
In view of the Monna-Springer theorem 8.4 \cite{roo}
there exists the non-trivial Haar measure $\mu : Bco({\bf L})\to \bf Q_p$
such that $\mu (B({\bf L},0,1))=1$,
where $Bco({\bf L})$ is the algebra of clopen (closed and open)
subsets in $\bf L$. By the definition of the Haar measure
$\mu (z+A)=\mu (A)$ for each clopen compact subset $A$ in $\bf L$ and each
$z\in \bf L$ in the non-Archimedean case and for each Borel subset
$A$ and each $z\in \bf L$ in the classical case 
(that is, for $X$ over fields $\bf R$ or $\bf C$).
The measure $\mu $ has an extension from $Bco({\bf L})$ on the 
$\sigma $-field $Bf({\bf L})$ of Borel subsets. We say that $A$ and 
$B \in 2^X_{\delta }$ are $\Upsilon $-equivalent and write
$A\Upsilon B$ if and only if 
$\int_{\pi (A\bigtriangleup B)}f(x)\mu (dx)=0$ for each 
$\pi \in \sf P$ and $\mu $-measurable $f$.
\par {\bf Theorem.} {\it Let $X$ be a $\kappa $-normed space.
Then there exists a one-to-one continuous mapping
$\Theta $  from $S_X/\Xi $ into either ${\bf R}^{d(X)}$, when 
${\bf K}=\bf R$ or ${\bf K}=\bf C$, or $\Psi $ from $S_X/ \Upsilon $
into ${\bf Q_p}^{d(X)}$
respectively, where $d(X)$ is the topological density of $X$,
$\bf p$ is the prime number such that the field of $\bf p$-adic numbers
${\bf Q_p}$ is not contained 
in $\bf K$. Moreover, the restrictions $\Xi |_{2^X_o}$ and
$\Upsilon |_{2^X_o}$ coincide with the equality of sets.}
\par {\bf Proof.} For ${\bf K}=\bf C$ 
we take ${\bf L}=\bf R$ and ${\bf L}=\bf R$ for ${\bf K}=\bf R$.
We consider the space $C^0({\bf L},Z)$ of continuous functions 
$f: {\bf L}\to Z$, where either $Z=\bf R$ or $Z=\bf Q_p$ respectively.
We take a countable family $\sf U$ in $C^0({\bf L},Z)$
satisfying two conditions:
\par $(i)$ the restriction ${\sf U}|_V$ is dense in
$C^0(V,Z)$ for each canonical closed compact subset $V$ in $\bf L$,
\par $(ii)$ either $\| \mu ^f \| <\infty $ in the case ${\bf L}=\bf R$
or $\sum_{j=1}^{\infty } \| \mu ^f|_{B({\bf L},x_j,1)} \| <\infty $
in the non-Archimedean case, where $\{ B({\bf K},x_j,1): j\in {\bf N} \} $
is a disjoint clopen covering of $\bf K$, $x_j\in \bf K$,
either $\| \mu ^f \| :=\int_{\bf R}|f(x)|\mu (dx)$ in 
the case ${\bf L}=\bf R$ or $\| \mu ^f |_A\| := 
\sup_{B\subset A, B\in Bco({\bf L})}
| \mu ^f(B) |$ for $A\in Bco({\bf L})$, $\mu ^f(A)=\int_Af(x)\mu (dx)$.
To continue the proof we need the following lemma.
\par {\bf 4.2. Lemma.} {\it Let $X$ be a $\kappa $-normable space.
If $\pi \in X^*_{\bf L}$ and 
$E\in 2^X_{\delta }$, then $\pi (E)\in Bf({\bf L})$, where 
$Bf({\bf L})$ is the $\sigma $-field of Borel subsets in $\bf L$.}
\par {\bf Proof.} In view of \S 2.4 \cite{eng} about quotient mappings
if $A$ is open in $X_{\bf L}$, then $\pi (A)$ is open in $\bf L$,
since $\bf L$ is spherically complete and
$\pi (A)=(A+Y)/Y$ due to the Hanh-Banach theorem
\cite{nari}, $Y:=coker (\pi )$ for each $\pi \in X^*_{\bf L}$.
In view of Lemma 2 \cite{lusmz} each $E\in 2^X_{\delta }$ is a 
countable intersection of open subsets $U_j$ such that 
$U_j\supset U_{j+1}$ for each $j\in \bf N$.
Then $\pi (E)=\bigcap_{j\in \bf N}((U+j+Y)/Y),$
since $\pi (E)\subset \bigcap_j \pi (U_j)$,
$[z]\in (E+Y)/Y$ if and only if there exists $z\in  E$ such that 
$(z+Y)/Y=[z]$ if and only if $z\in \bigcap_jU_j$ and 
$(z+Y)/Y=[z]\in (U_j+Y)/Y$ for each $j$. Therefore, $\pi (E)$ is the 
Borel subset in $\bf L$.
\par Continuation of the {\bf proof} of Theorem 4.1. In view of Lemma 4.2 
for each $E\ne C\in 2^X_o$ there exists $\pi \in \sf P$
such that $\pi (E)\ne \pi (C)$ and $f\in \sf U$ for which
$\mu ^f(\pi (E))\ne \mu ^f(\pi (C))$. 
Then $\mu ^f(\pi (E))=0$ for each $\pi \in \sf P$
and each $f\in \sf U$ is equivalent to $\mu (\pi (E))=0$
for each $\pi \in \sf P$ in the classical case, 
since $\mu ^f$ is absolutely continuous
relative to $\mu $ for each $f\in C^0({\bf L},Z)$.
By the diagonal mapping theorem 2.3.20 \cite{eng}
and Lemma 2 \cite{lusmz} we have a one-to-one continuous mapping 
$$F:=\Delta _{(\pi \in {\sf P}, f\in {\sf U})} 
\mu ^f(\pi (*)): W\to \prod_{(\pi \in {\sf P}, 
f \in {\sf U})}Y_{(\pi ,f)} ,$$ 
where either all $Y_{(\pi ,f)}$ 
are equal to $\bf R$ for $W=(2^X_{\delta }/\Xi )$
or all equal to $\bf Q_p$ for $W=(2^X_{\delta }/\Upsilon )$.
The minimal cardinality of $\sf P$ which separates points in 
$X_{\bf L}$ is equal to the topological density $d(X_{\bf L})$ 
of the space $X_{\bf L}.$ Evidently $d(X)=d(X_{\bf L})$,
since $X$ and $X_{\bf L}$ are homeomorphic as 
topological spaces (without linear structure).
Since $card ({\sf U})=\aleph _0$ and $d(X)\ge \aleph _0$,
then $card ({\sf U}) card (X)= card (X)$.
From the proof it follows, that for $Z={\bf L}=\bf R$ 
we can take instead of
$\sf U$ satisfying Conditions $4.1.(i,ii)$ simply ${\sf U}=\{ g \} $,
where $g(x)=e^{-x^2}$ for each $x\in \bf R$.
\par {\bf 4.3. Corollary.}{\it There exists a one-to-one
continuous mapping $\Theta $ of $2^X_o$ and $2^X_{\delta }/\Xi $
into $[0,1]^{d(X)}$ for ${\bf K}=\bf R$ 
and ${\bf K}=\bf C$ or $\Psi $ of $2^X_o$ and $2^X_{\delta }/\Upsilon $
into ${\bf Z_p}^{d(X)}$ in the non-Archimedean 
case.}
\par {\bf Proof.} Take in the proof of Theorem 4.1 in Condition 
$(ii)$ $\| \mu ^f \| \le 1$ for the real or complex field
and $\sum_j \| \mu ^f |_{B({\bf L},x_j,1)} | \le 1 $
in the non-Archimedean case.
\par {\bf 4.4. Note.} As shows a particular case of a normed space $X$ 
such embeddings $\Theta $ and $\Psi $ produce 
in general the weaker topologies
inherited from the Tychonoff product topology in $W:=\prod_{(\pi 
\in {\sf P}, f\in {\sf U})} {\bf L}$, than the initial one in the metric 
space $(2^X_o, D)$, since $W$ is not metrizable for $d(X)>\aleph _0$. 
In general $F$ is only continuous, but $F^{-1}$ may be discontinuous.
\par The topology in $2^X_o$ or in $2^X_{\delta }/\Xi $ 
or in $2^X_{\delta }/\Upsilon $ inherited from 
$W$ we call the h-weak topology and denote $\tau _w$. 
\par In the non-Archimedean case the $\kappa $-norm $\rho $
induces the ultrametric $D$ in $2^X_o$ or $2^X_{\delta }$, 
hence $(2^X_o,D)$ and $(2^X_{\delta },D)$ are totally disconnected.
\par {\bf 4.5. Corollary.} {\it If $X$ is over a non-Archimedean field, 
then $2^X_o$ and $2^X_{\delta }/\Upsilon $ 
are totally disconnected in the h-weak topology.}
\par {\bf 4.6 Corollary.} {\it $(2^X_o,\tau _w)$ 
and $(2^X_{\delta }/\Xi ,\tau _w)$ and $(2^X_{\delta }/\Upsilon ,\tau _w)$
have compactifications contained either in $[0,1]^{d(X)}$ 
for ${\bf K}=\bf R$ and ${\bf K}=\bf C$ or in ${\bf Z_p}^{d(X)}$
for the non-Archimedean field $\bf K$.}
\par {\bf 4.7. Note.} In accordance with Theorem 2.3.23 \cite{eng} there
exists a homeomorphic embedding of $(S_X,D)$ into $I^{w(S_X)}$,
where $w(S_X)$ is the topological weight of $S_X,$ where 
either $S_X=2^X_o$ or $S_X=2^X_{\delta }$, $I=[0,1]$.
Therefore, $w(S_X)\ge d(X)$, since $w(I^{\sf m})=\sf m$ 
for each ${\sf m}\ge \aleph _0$. 
Certainly equivalence relations $\Xi $ and $\Upsilon $
are different, since in the classical case $\mu $ is atomless
and in the non-Archimedean case $\mu $ on $Bf({\bf L})$
is purely atomic (because of Ch. 7 \cite{roo}).
Example 4.8 below and 
the results given above show, that functionals $\{ \mu ^f(\pi 
(*)): f, \pi \} $ do not separate points $E\in 2^X_{\delta }$ 
and closed subsets $F\subset 2^X_{\delta }$ with $E\notin F$,
that is, there are such $E$ and $F$, which are not separated.
\par {\bf 4.8. Example.} $(2^{\bf R^n}_o,D)$ is not separable and
it is dense in $(2^{\bf R^n}_{\delta },D)$, since countable unions
of closed parallelepipeds $\{ \prod_{i=1}^n [a_i,b_i] \} $
with $a_i, b_i \in \bf Q$ are dense in $2^{\bf R^n}_o$.
It can be used the Souslin number $hc(Y)$ of a topological space 
$Y$, which is the least cardinal such that every subset $A$ of $Y$ 
consisting exclusively of isolated points has $card (A)\le hc (Y)$, 
so $w(Y)\ge hc(Y)$ (see \S 1.7.12 \cite{eng}). On the other hand,
$2^{\bf R^n}_{\delta }$ has the metric inherited from $\bf R^n$.
Therefore,
$$\aleph _0^{\aleph _0}={\sf c}=w(2^{\bf R^n}_{\delta })=w(2^{\bf R^n}_o)
>w({\bf R^n})=\aleph _0 .$$
In the non-Archimedean case of the locally compact infinite field
with the non-trivial valuation countable unions of balls
$ \{ B({\bf K^n},x_j,r_j) : j\in \Omega , r_j \in \Gamma _{\bf L} \} $
are dense  in $2^{\bf L^n}_{\delta }$, where $\Omega $ is a countable 
dense subset in $\bf L^n$ and $\Gamma _{\bf L}:= \{ |z|_{\bf L}: 
0\ne z\in {\bf L} \} $ is discrete in $(0,\infty )$.
Therefore, $w(2^{\bf L^n}_{\delta })=w(2^{\bf L^n}_o)=\sf  c,$
since there exists a family of the cardinality $\aleph _0$ of
disjoint balls in $\bf R^n$ and in $\bf L^n$ with 
$\inf_{j\in \Omega }r_j>0$, hence $c({\bf L^n})=\aleph _0$ and
$hc(2^{\bf L^n}_o)=\sf c$.
\par {\bf 4.9. Theorem.} {\it There exists a homeomorphic embedding
of $2^X_{\delta }$ either into $I^{2^{w(X)}}$ for a $\kappa $-normable
topological vector 
space $X$ over $\bf R$ or $\bf C$, or into $ \{ 0, 1 \} ^{2^{w(X)}}$ for 
$X$ over a non-Archimedean infinite field with a non-trivial valuation,
where $ \{ 0,1 \} $ is a two-element dicrete set.}
\par {\bf Proof.} Each open subset is a union of some subfamily of elements
of a base, hence using closures of open subsets
we get $card (2^X_o)\le w(X)^{w(X)}$. Therefore, $w(2^X_o)\le 
2^{w(X)^{w(X)}}=2^{2^{w(X)}}$, since $w(X)\ge \aleph _0$.
In view of Lemma 2 \cite{lusmz} $w(2^X_{\delta })\le w(2^X_o)^{\aleph _0}$.
For the $\kappa $-normed space $(X,2^X_{\delta },\rho )$
the base of the topology of $2^X_o$ has the cardinality no greater, than
$w(X)^{w(X)}=2^{w(X)}$, since $w(X)\ge \aleph _0$ and
$\rho (x,C)$ is uniformly continuous by $x\in X$ for each fixed
$C\in 2^X_{\delta }$.
\par {\bf 4.10. Remark.} Then $w(X)\le w(2^X_{\delta })\le 2^{w(X)},$
since $w(2^X_{\delta })\ge w(X)$ and $card (2^X_o)\ge w(X)$.
\section{Applications of $\kappa $-normed spaces.}
\par {\bf 5.1. Theorem.} {\it Let $X$ be a complete
$\kappa $-normable locally convex space
over the field ${\bf K}=\bf R$ or $\bf C$ and let be given
a continuous mapping
$f: {\bf R}\times (2^X_{\delta }\cup X)\to 2^X_{\delta }\cup X$ such that
$f(t,2^X_{\delta })\subset 2^X_{\delta }$, $f(t,X)\subset X$ and
$\rho _X(f(t,x), f(t,A))\le C_A\rho _X(x,A)$ for each $A\in 
2^X_{\delta }$ and $x\in X$, where a constant $C_A>0$ may depend on
$A$, $f(t,M)=M$ for each $t\in \bf R$ and $\bf K$-linear subspace
$M$ from the family $\sf F$ such that ${\sf F}\subset 2^X_{\delta }$,
$M_1\cap M_2\in \sf M$ for each $M_1$ and $M_2\in \sf F$,
$\bigcap_{M\in \sf F}M=\{ 0\}$. Then the differential equation
$dx(t)/dt=f(t,x(t))$ (or $dA(t)/dt=f(t,A(t))$)
with the initial condition $x(0)=x_0\in X$ 
(or $A(0)=A_0\in 2^X_{\delta }$
for a convex balanced subset $A_0-x_0$ for some $x_0\in X$)
has the unique solution $x: {\bf R}\to X$ (or $A: {\bf R}
\to 2^X_{\delta }$ for ${\sf F}=2^X_{\delta }$
correspondingly).}
\par {\bf Proof.} $(I)$. Let $M$ be a $\bf K$-linear
subspace such that $M\in \sf F$ (see the end of \S 2.3). Then
$\rho _X(x,M)=\| x+M\|_{X/M}$ is the norm in the quotient space
$X/M$, where $x+M=\theta _M(x)$, $\theta _M: X\to X/M$ is a quotient 
mapping,
$x\in X$, $(x+M)\in X/M$. The function $f_M(t,x)=f(t,x)+M$ 
satisfies the following condition:
$\| f_M(t,x)\| _{X/M}=\rho _X(f(t,x),f(t,M))\le C_M\rho _X(x,M)=
C_M\| x+M\|_{X/M}$, since $f(t,M)=M$. Moreover, 
$f_M: X/M\to X/M$ is continuous.
The equation $dx_M(t)/dt=f_M(t,x_M(t))$ with the initial condition
$x_M(0)=x_0+M$ has the unique solution
$x_M(t)$ in the Banach space $X/M$.
Evidently that in the particular case, when 
${\sf F}=2^X_{\delta ,l}$ is the family of all
$\bf K$-linear subspaces belonging to $2^X_{\delta }$,
the conditions imposed on $\sf F$ are satisfied, therefore
the conditions of this theorem are correct.
If $M_1$ and $M_2\in \sf F$, then $x_{M_1}(t)\cap x_{M_2}(t)=
x_{M_1\cap M_2}(t)$, since $(f(t,x)+M_1)
\cap (f(t,x)+M_2)=f(t,x)+(M_1\cap M_2)$.
We have $\bigcap_{M\in \sf F}x_M(t)=x(t)\in X$ and $\bigcap_{M\in \sf F}
f_M(t,x)=f(t,x)$, consequently, $dx(t)/dt=f(t,x)$.
\par $(II)$. Let $A\in 2^X_{\delta }$ be convex and balanced, 
then $M:=(\bigcap_n\lambda _nA)
\in 2^X_{\delta ,l}$, where $0\ne \lambda _n\in \bf K$ 
for each $n\in \bf N$ and $\lim_{n\to \infty }\lambda _n=0$. 
Consequently, $\theta _M(A)=(A+M)
\in 2^{X/M}_{\delta }$, in addition $A+M$ is bounded in the Banach space
$X/M$, since $X$ is complete. From the inequality 
$\rho _X(f(t,x),f(t,A))\le C_A\rho _X(x,A)$ 
it follows that $D(f(t,A_1), f(t,A_2))=\bar \rho _X(f(t,A_1), f(t,A_2))+$
$\bar \rho (f(t,A_2), f(t,A_1))$ $\le [C_{A_1}\bar \rho _X(A_2,A_1)+C_{A_2}
\bar \rho _X(A_1,A_2)]$ $\le D(A_1,A_2)C$, where $C=\max (C_{A_1},C_{A_2})$. 
Then $\int_{t_1}^{t_2}[f(\tau ,A)+M]d\tau =$ $\int_{t_1}^{t_2}f(\tau ,A)
d\tau +M$, since $tM=M$, where $t_2>t_1$, $t:=t_2-t_1$, $\int_{t_1}^{t_2}E(
\tau )d\tau $ $:=\{ \phi (t_1,t_2):=\int_{t_1}^{t_2} \psi (\tau )d\tau $,
$\psi (\tau )\in E(\tau ) \} $, $E(\tau ):=f(\tau ,E)=\{ f(
\tau ,x): x\in E \} $, $\psi (\tau )=f(\tau ,x)$.
In Lemma 8 \cite{lusmz} it was proved that the function
$D(A,B):={\bar \rho }(A,B)+{\bar \rho }(B,A)$ on $(2^X_{\delta })^2$
is the metric satisfying additional conditions
$D(A{\hat +}B,C{\hat +}E)\le D(A,C)+D(B,E)$ (or with the maximum
instead of the sum in the non-Archimedean case) $D(\lambda A,\lambda C)=
|\lambda |D(A,C)$ for each $0\ne \lambda \in \bf K$, where
$A{\hat +}B:=cl(A+B)$.
Using the approximation of continuous functions by simple (step) functions 
in the space $L^1([a,b],X)$ for $a\le t_1<t_2\le b$ we get:
\par $D(\int_{t_1}^{t_2}f_M(\tau ,A_1)d\tau , 
\int_{t_1}^{t_2}f_M(\tau ,A_2)d\tau )$ $\le C(t)D(A_1,A_2)$, \\
that is, the mapping $F_M(t_1,t_2,A):=\int_{t_1}^{t_2}f_M(\tau ,A)d\tau $
is contracting for $Ct<1$. Then in each interval $[t_1,t_2]
\ni 0$ with $t=t_2-t_1<C^{-1}$ there exists the unique solution
$(A+M)(t)$ in $X/M$, where $(A+M)(t)=A(t)+M$. 
Since $M\subset A$, then $A+M=A$, that is, $A(t)$ is a solution
with the initial condition $A_0=A(0)$, where $M=\bigcap_n\lambda _nA_0$.
Using finite coverings of segments
$[a,b]$ and solving analogous tasks with initial conditions
in intermediate points and sewing solutions we get 
the solution in each segment $[a,b]\ni 0$.
Evidently, the second part of this theorem is accomplished 
for more general functions
$f$ with $f(t,M)=M$ for each $t\in \bf R$ and a given one subspace
$M$ such that $M\in 2^X_{\delta ,l}$, where $M$ is completely defined
by the initial condition.
\section{Approximations of representations of the interval orders 
with the help of $\kappa $-normed spaces.}
\par {\bf 6.1. Remarks.} Mappings of interval orders play 
an important role in mathematical economics.
In papers \cite{bois1,bois2,bois3,is1} a problem of representations of 
interval orders was began to study naively 
with the help of pairs of semicontinuous 
and continuous functions. It was proved the existence of such pairs and 
examples were given, but this field remains investigated only a little.
The problem of searching of all such pairs
and investigations of cases, when
it is possible to use one function instead of the pairs remains actual.
For economic theories it is important to search for the best
approximation. 
\par If $T$ is a preordered set or a set with an interval 
order, then its representation is searched by a function $f: T
\to \bf R$, which preserve order, or a pair of functions
$f, g :T\to \bf R$ such that $g(x)\ge 0$
for each $x\in T$ and $x_1<x_2$ if and only if 
$v(x_1)+\sigma (x_1)<v(x_2)$. This section is devoted to
the investigation of this problem with the help of
$\kappa $-normed spaces.
\par Recall that by the definition the mapping $f: T\to \bf R$ 
represents $T$ if and only if $f: T\to \bf R$ is order-preserving.
Certainly not for every $T$ such $f$ exists. Therefore,
approximations are frequently more valuable, than the precise results,
which may be nonexistent.
\par {\bf 6.2. Definitions and notes.} Let $T$ be a preordered set
and ${\sf U}_c$ be a family of linearly ordered countable subsets $t$ in $T$
and ${\sf B}_c$ be a subfamily of $t\in \sf U$ such that there are
$a_t$ and $b_t\in T$ with $a_t\le x$ and $x\le b_t$
for each $x\in t$. When the condition of countability of $t$ is 
dropped we denote the corresponding families
by $\sf U$ and $\sf B$. The family of closed subspaces $t$ in $T$ 
we denote by $\sf C$. There can be considered another families
$\sf F$ of subsets in $T$.
\par Let $C_b(t,{\bf R})$ denotes a family of functions $f: t\to \bf R$ 
such that $\| f \| :=\sup_{x\in t} |f(x)|<\infty $.
If $T$ is supplied with a topology $\tau $ and each $t\subset T$ is 
considered in a topology inherited from $T$, then 
$C^0_b(t,{\bf R})$ denotes a subspace
of continous functions $f$ in $C_b(t,{\bf R})$.
\par {\bf 6.3. Theorem.} {\it For each family $\sf F$ there
are $\kappa $-normable spaces $(X_{\sf F},S_{X_{\sf F}})$ 
both for the cases
of a family of all canonical closed subsets $S_{X_{\sf F}}=2^{X_{\sf F}}_o$
and a family of closed $G_{\delta }$-subsets 
$S_{X_{\sf F}}=2^{X_{\sf F}}_{\delta }$.}
\par {\bf Proof.} Each space $C_b(t,{\bf R})$ and $C^0_b(t,{\bf R})$ is 
normed and hence $\kappa $-normed. The first space is the particular case
of the second, when the set $t$ is supplied with the discrete topology.
Then $P_{\sf F}:=\prod_{t\in \sf F}Y_t$ is the $\kappa $-normable space 
by a $\kappa $-norm $\rho $ on $(X_{\sf F},2^{X_{\sf F}}_o)$,
where $P_{\sf F}$ is in the product locally convex topology
(see Theorem 14 \cite{lusmz}). If we take $S_{\sf F}$
as a subspace of $P_{\sf F}$ satisfying Conditions
$(a-d)$ of Proposition 13 \cite{lusmz}, then we get
a $\kappa $-norm $\rho $ on $(S_{\sf F},2^{S_{\sf F}}_{\delta })$.
The topology of the latter locally convex space 
$S_{\sf F}$ satisfies conditions of Theorem 3.2 above.
\par In view of Theorem 6 \cite{lusmz} the quotient space of 
the $\kappa $-normed space is $\kappa $-normed, then it can be taken 
into account an equivalence relation $\sf E$ in $P_{\sf F}$ 
and $S_{\sf F}$ generated by equalities 
$$(i)\mbox{ } Y_{t_1\cap t_2}=Y_{t_1}\cap Y_{t_2}$$
in the case of $C_b(t,{\bf R})$. On the other hand,
in the case of $C_b^0(t,{\bf R})$ the equality $(i)$ is satisfied, when
each $f\in C_b^0(t_1\cap t_2,{\bf R})$ has a continuous bounded
extension on $t_1\cup t_2$. This is the case, for example,
for $\sf F$ consisting of closed subspaces $t$ in a normal 
topological space $T$ (see Tietze-Uryson theorem 2.1.8 \cite{eng}).
When Condition $(i)$ is satisfied we can take either 
$X_{\sf F}=P_{\sf F}/{\sf E}$ or $X_{\sf F}=S_{\sf F}/{\sf E}$, 
which fits conditions of this theorem.
\par {\bf 6.4. Theorem.} {\it Let $(X,S_X,\rho )$ be a $\kappa $-normed 
space and $D$ be a metric in $S_X$ induced by $\rho $. Then for each $C_v\in S_X$
and each $\epsilon >0$ there exists $\emptyset \ne C_{\sigma }
\in S_X$ such that $D(cl(C_v+C_{\sigma }),C_v)<\epsilon $.}
\par {\bf Proof.} In view of Corollary 5 \cite{lusmz}
for each $C_v\in S_X$ and each $\epsilon >0$ there exists
an open neighbourhood $U$ of zero in $X$ such that
$0\le \rho (f,C_v)-\rho (f,cl(C_v+U))<\epsilon $
for each $f\in X$. We take $C_{\sigma }\in S_X$ such that $\emptyset \ne 
C_{\sigma }\subset U$  and $0\in C_{\sigma }$. Then  $\rho (f,
cl(C_v+C_{\sigma }))\ge \rho (f,cl(C_v+U))$ for each $f\in X$ due to 
the monotonicity axiom of $\rho $. Therefore, $0\le \rho (f,C_v)-\rho (f,
cl(C_v+C_{\sigma }))<\epsilon $ for each $f\in X$. If 
$f\in cl(C_v+C_{\sigma })$, then $\rho (f,cl(C_v+C_{\sigma }))=0$
due to the inclusion axiom. On the other hand,
$\rho (f,C_v)\ge \rho (f,cl(C_v+C_{\sigma }))$
due to the monotonicity axiom. Hence ${\bar \rho }(C_v,cl(C_v+C_{\sigma }))
=\sup_{f\in C_v}\rho (f,cl(C_v+C_{\sigma }))=0$. In view of the triangle 
inequality axiom $\rho (f,C_v)\le \rho (f,cl(C_v+C_{\sigma }))+
{\bar \rho }(cl(C_v+C_{\sigma }),C_v).$ Then $0\le \rho (f,C_v)<\epsilon $ 
for each $f\in cl(C_v+C_{\sigma })$, consequently, 
$D(C_v,cl(C_v+C_{\sigma }))={\bar \rho }(C_v,cl(C_v+C_{\sigma }))+
{\bar \rho }(cl(C_v+C_{\sigma }),C_v)=\sup_{f\in cl(C_v
+C_{\sigma })}\rho (f,C_v)\le \epsilon .$
\par {\bf 6.5. Theorem.} {\it Let $T$ be a preordered set
or $T$ be a set with an interval order and $\sf F$ be a family of subsets
of $T$ such that there exists a $\kappa $-normed space
$(X_{\sf F},2^{X_{\sf F}}_{\delta },\rho )$. Suppose there is 
a function $f: T\to \bf R$
which represents each $t\in \Lambda $ and $f|_t$ is bounded, 
where $card (\Lambda )<\aleph _0$ and $\Lambda \subset \sf F$.
\par  $(1)$. Then there exists $C\in 2^{X_{\sf F}}_{\delta }$ 
such that $f\in C$. 
\par $(2)$. Moreover, $C$ can be chosen such that 
$\pi _t(C)=f|_t$ for each $t\in \Lambda $.}
\par {\bf Proof.} Let $X_{\sf F}$ be a $\kappa $-normed space
from Theorem 6.3 with $S_{X_{\sf F}}=2^{X_{\sf F}}_{\delta }$. 
When $\Lambda $ is a finite family
we can take $2^{X_{\sf F}}_o$ also instead of $2^{X_{\sf F}}_{\delta }$,
but the Condition $(2)$ in general may be unsatisfied.
Therefore, we have the $\kappa $-normed space
$(X_{\sf F},2^{X_{\sf F}}_{\delta },\rho )$.
In view of Theorem 6 \cite{lusmz} the $\kappa $-norm
$\rho $ on $(X_{\sf F},2^{X_{\sf F}}_{\delta })$ generates a 
norm on $Y_t$ by the formula ${\tilde p}(\pi _t(f))=\rho (f,Y_t)$,
since $Y_t\in 2^{X_{\sf F}}_{\delta }$, where $\pi _t: X_{\sf F}\to Y_t$ is 
the quotient mapping. This norm is equivalent to the initial one.
On the other hand, $\pi _t(f)=f|_t$ and $\{ \pi _t(f) \} $
is a $G_{\delta }$-subset in $Y_t$. Hence $C=\bigcap_{t\in \Lambda }
\pi _t^{-1}( \{ \pi _t (f) \} )$ is in $S_{X_{\sf F}}$ and $f\in C$.
This $C$ also satisfies Condition $(2)$ for $S_{X_{\sf F}}
=2^{X_{\sf F}}_{\delta }$.
\par In the case of $card (\Lambda )<\aleph _0$ and $S_{X_{\sf F}}
=2^{X_{\sf F}}_o$ we can take balls
$B_t:= \{ y: y\in Y_t, p_t(y_t-\pi _t(f))\le \epsilon _t \} $
with $\infty >\epsilon _t>0$ for each $t\in \Lambda $,
then $C=\bigcap_{t\in \Lambda }\pi _t^{-1}(B_t)$ satisfies
Condition $(1)$. For $card (\Lambda )=\aleph _0$ and
$S_{X_{\sf F}}=2^{X_{\sf F}}_{\delta }$ we can take
$C=\bigcap_{t\in \Lambda }\pi _t^{-1}(B_t)$ satisfying Condition $(1)$.
\par {\bf 6.6. Corollary.} {\it Let ${\sf F}\in \{ {\sf U},{\sf U_c},
{\sf B}, {\sf B_c},{\sf C} \} $, then there exists a $\kappa $-normed 
space $(X_{\sf F},2^{X_{\sf F}}_{\delta },\rho )$ such that
for each representation function $f$ of a preordered or an 
interval ordered set $T$ on a countable subfamily $\Lambda $
of a family ${\sf F}\in \{ {\sf U}, {\sf U_c}, {\sf B}, {\sf B_c}, 
{\sf C} \} $, for which the restriction $f|_t$ is bounded and  
represents $t$ for each $t\in \Lambda $, there exists
$C\in 2^{X_{\sf F}}_{\delta }$ with $f\in C$
and $\pi _t(C)=f|_t$ for each $t\in \Lambda $.}
\par {\bf 6.7. Note.} For the family $\sf B_c$ the following
condtion: "$f|_t$ is bounded" is automatically satisfied, since
$-\infty <f(a_t)\le f(x)\le f(b_t)<\infty $ for each  $x\in t$,
where $a_t=\inf_{a\in  t}a$ and $b_t:=\sup_{b\in t}b$.
The family ${\sf U_c}$ corresponds to approximation of $f$ on 
each countable family of linearly ordered countable subsets.
\par The family $\sf C$ corresponds to approximation of $f$ on countable 
families of closed subsets. In particular there can be taken
a subfamily of compact subsets. This construction can be refined
(see below). In general in the segment $[a,b]\subset \bf R$
there can be embedded ${\sf c}=card ({\bf R})$ distinct subsets,
which are linearly ordered. For example, the ring $\bf Z_p$ of integer 
$p$-adic numbers $x$, that is, $x=\sum_{j=0}^{\infty }x_jp^j$,
where $x_j\in \{ 0,1,...,p-1 \} $, $p$ is the prime number.
The ring $\bf Z_p$ is linearly ordered by the relation
$x\triangle y$ if and only if $x_0=y_0$, $x_1=y_1,$...,
$x_n=y_n$ and $x_{n+1}<y_{n+1}$. As it is well-known
$\bf Z_p$ is totally disconnected and homeomorhic with
the Cantor set $\{ 0, 1 \} ^{\aleph _0}$, where $\{ 0, 1 \} $ is 
the discrete two-element set. The non-Archimedean metric in 
$\bf Z_p$ can correspond to economic models with the reciprocal dependence
on a parameter $x\in T$, $x\mapsto x^{-1}$. There are 
$\sf c$ pairwise distinct subsets in $[a,b]$, which are
homeomorhic with $\bf Z_p$. In $[a,b]$ there exist also
$\sf c$ pairwise distinct subsets (either countable
or) of the cardinality $\sf c$ with the linear order
inherited from the field $\bf R$.
\par {\bf 6.8. Theorem.} {\it Sets $C$ in Theorem 6.5 can be 
chosen such that 
\par $(i)$ for each $g\in C$ and $t\in \Lambda $ 
a restriction $g|_t$ is a nondecreasing function.}
\par {\bf Proof.} Let $C$ be a set from Theorem 6.5. We construct 
from it a new set denoted by $V$ satisfying Condition $(i)$.
Since $X_{\sf F}$ is a complete locally convex space and there exists
a topological dual space $X^*_{\sf F}$ of continuous linear functionals 
on  $X_{\sf F}$ such that $X^*_{\sf F}$ separates points of $X_{\sf F}$ 
by the Hahn-Banach theorem \cite{nari}. In particular
a mapping $G_{x_1,x_2}(f):=f(x_1)-f(x_2)$
for a fixed pair $\{ x_1,x_2 \} \subset T$ is a continuous 
linear functional on $X_{\sf F}$.
Then $C^+_{x_1,x_2}:= \{ g: g\in C, g(x_1) \le g(x_2) \} $
is a closed subset in $C$, since $g(x_1)\le g(x_2)$ is equivalent to
$G_{x_1,x_2}(g) \le 0$. On the other hand,
$C^{\epsilon }_{x_1,x_2}:= \{ g: g\in C, G_{x_1,x_2}(g)<\epsilon \} $
is open in $C$ for each $\epsilon >0$. Then $\bigcap_{n=1}^{\infty }
C^{1/n}_{x_1,x_2}=C^+_{x_1,x_2},$ consequently, $C^+_{x_1,x_2}\in 
2^{X_{\sf F}}_{\delta }$. We take $V=\bigcap_{(x_1<x_2; x_1, x_2
\in t, t\in \Lambda )} C^+_{x_1,x_2}$. Since $card (\Lambda )
\le \aleph _0$ and $card (\bigcup_{n=2}^{\infty }\Lambda ^n)=\aleph _0$, 
then $V\in 2^{X_{\sf F}}_{\delta }$. It is only necessary to verify
that $V\ne \emptyset .$ It is evident, when $f$ is given, since
$f\in C^+_{x_1,x_2}$ for each $x_1<x_2$ in $t$.
\par {\bf 6.9. Theorem.} {\it Let $\Lambda $ be a countable subfamily 
of $\sf F$, where ${\sf F}\in \sf U_c$. Suppose $\emptyset \ne C\in 
2^{X_{\sf F}}_{\delta }$ and $\pi _t(C)$ has a nonempty interior for each
$t\in \Lambda $, where $\pi _t$ are quotient mappings from Theorem 6.3,
then $\bigcap_{x_1<x_2; x_1, x_2\in t; t\in \Lambda }C^+_{x_1,x_2}=V
\ne \emptyset $ and $V\in 2^{X_{\sf F}}_{\delta }$.}
\par {\bf Proof.} In \S 6.8 it was proved that $V\in 
2^{X_{\sf F}}_{\delta }$. If $x_1<x_2$ in $t$, where $t\in \Lambda $, 
then $ \| G_{x_1,x_2} \| _{Y^*_t}\le 2$, hence $C^+_{x_1,x_2}$ is 
contained in $\pi _t^{-1}(U^0)$, where $U_t^0=U^0$ is the absolute polar 
of the ball $U=B(Y^*_t,0,2)$ in $Y^*_t$, $U^0\subset Y_t$,
$(Y_t,Y^*_t)$ is the dual pair. In view of the Alaoglu-Bourbaki theorem
(9.3.3) \cite{nari} $U^0$ is $\sigma (Y_t,Y^*_t)$-compact.
Since $\pi _t(C)$ contains the  ball $B(Y_t,g_t,r_t)$ we can take 
${\tilde C}\in 2^{X_{\sf F}}_{\delta }$ such that $\pi _t({\tilde C})
=B(Y_t,g_t,r_t)$ for each $t\in \Lambda $. In view of the Hahn-Banach 
theorem in its geometric form (8.5.3) \cite{nari} the ball $B(Y_t,g_t,r_t)$ 
is weakly closed, hence $\pi _t({\tilde C})$ is a compact subset of $U^0$.
In view  of Theorem 3.1 we can choose a family of seminorms
$p_t$ defining topology of $X_{\sf F}$ such that 
$\sup_{t\in \Lambda }p_t(x)=:p(x)<\infty $ is a continuous seminorm
in $X_{\sf F}$, since $card (\Lambda )\le \aleph _0$.
So instead of $t\in \Lambda $ we can take $t_0=\bigcup_{t\in \Lambda }t$. 
Then $\pi _{t_0}({\tilde C})$ is a weakly compact subset of $U^0_{t_0}$, 
where $U^0_{t_0}$ is for the space $Y_{t_0}$. For each $x_1<x_2<...<x_n$ 
there evidently exists $g: T\to \bf R$ such that $g\in X_{\sf F}$ and 
$g(x_1)<g(x_2)<...<g(x_n)$. Then $\pi _{t_0}({\tilde C}^+_{x_1,x_2})$ 
are weakly closed subsets of $\pi _{t_0}({\tilde C})$. Since
$(\bigcap_{x_1<x_2; x_1, x_2\in t; t\in \Lambda }\pi _{t_0}(
{\tilde C}^+_{x_1,x_2}))\ne \emptyset $, then $\bigcap_{(x_1<x_2; x_1, 
x_2 \in t; t\in \Lambda )}{\tilde C}^+_{x_1,x_2}\ne \emptyset $.
\par {\bf 6.10. Note.} The construction above can be generalized to 
encompass the case of unbounded functions.
For this we substitute $ \| f \| =\sup_{x\in t}|f(x)|$ on
another seminorm $ \| f \| =\sup_{x\in t}|f(x)|\phi (x)$,
where $\phi : t\to (0,\infty )$ is a fixed function such that
$\phi (x)>0$ for each $x\in t$. We can take different $\phi $ and 
$t\in \sf F$ in particular such that $1/\phi $ is unbounded. 
With sufficiently large family $\{ \phi \} $ we can take into account 
unbounded $f: T\to \bf R$.
\par When $T$ is either a subset of $\bf R$ in a topology inherited 
from $\bf R$ or $T\subset \bf  Q_p$ in a topology iherited from $\bf Q_p$, 
then we can consider also functions satisfying the Lipshitz condition
$$(i)\mbox{ }|g(x_1)-g(x_2)|\le C_1|x_1-x_2|,$$ where $C_1=const>0$
and also satifying condition
$$(ii) \mbox{ }(g(x_2)-g(x_1))\ge C_2(x_2-x_1)$$
for each $x_2>x_1$ in $t$, where $t\in \Lambda $
(in $(ii)$ for $T\subset \bf R$). Then Theorem 6.9 can be modified 
that to satisfy Condition $(ii)$, where $C_2$ can be dependent on $t$.

\par Theoretical Department, Institute of General Physics of
Russian Academy of Sciences, 
\par str. Vavilov 38, Moscow, 117942, Russia
\end{document}